\documentclass[12pt,dvips]{amsart}
\usepackage{amsfonts, amssymb}
\usepackage[usenames, dvipsnames]{color}
\usepackage{diagrams}
\usepackage{hyperref}

\newtheorem{lemma}{Lemma}[section]

\newtheorem{corollary}{Corollary}[section]

\newtheorem{proposition}{Proposition}[section]

\newtheorem*{declaration}{Declaration}
\newtheorem*{acknowledgements}{Acknowledgements}

\newsavebox{\tempbox}
\newlength{\templength}
\newcommand{\vcenterbox}[1]
   {\sbox{\tempbox}{#1}
    \settowidth{\templength}{\usebox{\tempbox}}
    \parbox{\templength}{\usebox{\tempbox}}}

\newdimen{\cellsize}
\newcommand\medboxes{\setlength{\cellsize}{14pt}\def\boxformat{}}

\medboxes
\newsavebox{\cellcontent}
\def\hidehrule#1#2{\kern-#1
  \hrule height#1 depth#2 \kern-#2 }%
\def\hidevrule#1#2{\kern-#1{\dimen\cellcontent=#1%
    \advance\dimen\cellcontent by#2\vrule width\dimen\cellcontent}\kern-#2 }%
\def\makeblankbox#1#2{\hbox{\lower\dp\cellcontent\vbox{\hidehrule{#1}{#2}%
    \kern-#1 
    \hbox to \wd\cellcontent{\hidevrule{#1}{#2}%
      \raise\ht\cellcontent\vbox to #1{}
      \lower\dp\cellcontent\vtop to #1{}
      \hfil\hidevrule{#2}{#1}}%
    \kern-#1\hidehrule{#2}{#1}}}}
\newcommand\cellify[1]{\defaultcell%
\sbox{\cellcontent}{\vbox to \cellsize{%
\vfill%
\hbox to \cellsize{\hfill$\boxformat #1$\hfill}
\vfill}}%
\rlap{\drawnbox}
\usebox{\cellcontent}}
\newcommand\tableau[1]{\vtop{\let\\\cr
\baselineskip -16000pt \lineskiplimit 16000pt \lineskip 0pt
\ialign{&\cellify{##}\cr#1\crcr}}}
\newcommand\defaultcell{\gdef\drawnbox{
\makeblankbox{0.2pt}{0.2pt}
}}
\newcommand\graycell{\gdef\drawnbox{%
\rlap{\color{Gray}\vrule width \cellsize height \cellsize}%
\makeblankbox{0.2pt}{0.2pt}
}}
\newcommand\thickcell{\gdef\drawnbox{
\makeblankbox{0.2pt}{0.1\cellsize}%
}}
\newcommand\vdotscell{\gdef\drawnbox{\kern-1.6pt\vbox{\baselineskip=4pt\lineskiplimit=0pt\hbox{}\hbox{.}\hbox{.}\hbox{.}}}}
\begin{document}

\title{Symmetric Functions and Macdonald polynomials}
\author{Robin Langer}	
\date{September 15 2008}
\maketitle
\begin{center}
Submitted in total fulfilment of the requirements of the degree of Master of Mathematics by Research at the University of Melbourne
\end{center}

\newpage
\begin{abstract}
The ring of symmetric functions $\Lambda$, with natural basis given by the Schur functions, arise in many different areas of mathematics. For example, as the cohomology ring of the grassmanian, and as the representation ring of the symmetric group. One may define a coproduct on $\Lambda$ by the plethystic addition on alphabets. In this way the ring of symmetric functions becomes a Hopf algebra. The Littlewood--Richardson numbers may be viewed as the structure constants for the co-product in the Schur basis. The first part of this thesis, inspired by the umbral calculus of Gian-Carlo Rota, is a study of the co-algebra maps of $\Lambda$. The Macdonald polynomials are a somewhat mysterious qt-deformation of the Schur functions. The second part of this thesis contains a proof a generating function identity for the Macdonald polynomials which was originally conjectured by Kawanaka.
\end{abstract}

\newpage

\begin{declaration}
The first part of this thesis is entirely the author's own work. The idea for the second part was suggested by Ole Warnaar who also carried out the derivation of equation \ref{ident} on pages 59--60 and verified proposition \ref{final}. The proof of lemma \ref{inflim} was suggested by Paul Zinn-Justin. Steps two and three in the proof of the second part of the thesis are the author's own work.
\end{declaration}

\newpage

\begin{acknowledgements}
The author would like to thank Ole Warnaar for suggesting the subject of Macdonald polynomials, Peter Forrester, Alain Lascoux, Paul Zinn-Justin and the two anonymous referees for useful comments and suggestions on the manuscript as well as the University of Montreal for hosting an interesting lecture series and workshop on Macdonald Polynomials and Combinatorial Hopf Algebras.
\end{acknowledgements}

\newpage

\hfill \break

\hfill \break

\hfill \break

\hfill \break

\begin{center}
{\it for Mrs Mclaren}
\end{center}

\newpage

\setcounter{tocdepth}{3}{\small\tableofcontents}

\section{Symmetric functions of Littlewood--Richardson type}

The Schur polynomials $\{s_\lambda(X)\}$ arise in both geometry and representation theory as a natural basis for the ring of symmetric polynomials which we denote by $\Lambda$. The Littlewood--Richardson numbers $\{c^{\lambda}_{\mu \nu}\}$ are defined by:
$$ s_\mu(X)s_\nu(X) = \sum_\lambda c^{\lambda}_{\mu \nu} s_\lambda(X) $$

The plethystic addition of alphabets gives a coproduct structure on $\Lambda$, which induces, via the Hall inner product, the natural product structure on the dual space $\Lambda^*$ - which may be thought of as the ring of symmetric ``formal power series'' in some dual alphabet $Y$.

Since the Schur polynomials are self-dual with respect to the Hall inner product, the Littlewood--Richardson numbers may also be defined by:
$$ s_\lambda(X + X') = \sum_{\mu, \nu} c^{\lambda}_{\mu \nu} s_\mu(X) s_\nu(X') $$

A sequence of binomial type is a basis $p_n(x)$ for the one variable polynomial ring $\mathbb{Q}(x)$ with the property that:
$$ p_n(x+x') = \sum_{k=0}^n \binom{n}{k} p_k(x) p_{n-k}(x') $$
Sequences of binomial type arise as images, under co-algebra maps, of the standard basis $\{x^n\}$. Furthermore, each such co-algebra map, or {\it umbral operator} is associated to an a formal power series $f(z)$ which is invertible with respect to function composition.

The homomorphism from the one variable polynomial ring $\mathbb{Q}[x]$ to the ring of symmetric functions $\Lambda$ given by:
$$ \frac{x^n}{n!} \mapsto h_n(X) $$
preserves the coproduct structure. Similarly, the homomorphism from the ring of formal power series in one variable $\mathbb{Q}[[y]]$ to $\Lambda^*$ given by:
$$ y^n \mapsto m_{(n)}(Y) $$
preserves the product structure. 

In the special case where the alphabet $X$ contains a single variable, both these maps are in fact isomorphisms, and the Hall inner product becomes the bracket used in the umbral calculus of Gian-Carlo Rota:
$$ \langle f(y), p(x) \rangle = \mathcal{L}f(\mathcal{D})[p(x)] $$
Here $\mathcal{L}$ denotes the constant term operator, while $\mathcal{D}$ is the differential operator. 

The Schur functions may be expressed as either a determinant in the complete symmetric functions:
$$ s_\lambda(X) = \det( h_{\lambda_i + j - i}(X)) $$
or as a determinant in the elementary symmetric functions:
$$ s_\lambda(X) = \det( e_{\lambda'_i + j - i}(X)) $$

The main result of the first part of this thesis is the following. Let $\{p_n(x)\}$ be the sequence of binomial type arising as the image of $\{x^n\}$ under the umbral operator $U_f$ which is associated to the invertible formal power series $f(z)$. Let $\{r_n(X)\}$ be the image of $\{\frac{p_n(x)}{n!}\}$ under the embedding of $\mathbb{Q}[x]$ into $\Lambda$ just described. Then the vector space basis for $\Lambda$ defined by:
$$ P_\lambda(X) = \det( r_{\lambda_i + j - i}(X) ) $$
has the Littlewood--Richardson property:
$$ P_\lambda(X + X') = \sum_{\mu, \nu} c^{\lambda}_{\mu \nu} P_\mu(X) P_\nu(X') $$

Furthermore, suppose that $g(z)$ is the compositional inverse of $f(z)$ and that $\{q_n(x)\}$ is the sequence of binomial type arising as the image of the standard basis $\{x^n\}$ under the umbral operator $U_g$. Let $\{\tilde{r}_n(X)\}$ denote the image of $\{\frac{q_n(x)}{n!}\}$ under the embedding of $\mathbb{Q}[x]$ into $\Lambda$ and let:
$$c_n(X) = \omega[\tilde{r}_n(X)]$$
where $\omega$ is the natural involution on $\Lambda$ which maps the complete symmetric functions to the elementary symmetric functions. 

The aforementioned basis $\{P_\lambda(X)\}$ has an alternate expression of the form:
$$ P_\lambda(X) = \det( c_{\lambda_i + j - i}(X))$$
Every basis $\{P_\lambda(X)\}$ which has the Littlewood--Richardson property, and which also admits an expansion of the form:
$$ P_\lambda(X) = \sum_{\mu \subseteq \lambda} \Box_{\lambda \mu} s_\mu(X)$$
arises in this way.

In section 1.1 we review some classical results from the theory of symmetric functions. In section 1.2 we review classical results from the umbral calculus, in section 1.3 we prove the main theorem and in section 1.4 we give some examples.

\subsection{Symmetric Functions}

The main references for this section are Macdonald \cite{BLACK} and Lascoux \cite {LASCOUX}. Another good reference is Stanley \cite{STANLEY}. Also, see Bergeron's online notes (in French) \cite{BERG}. For connections with the geometry of grassmanians see Fulton \cite{BLUE} and Manivel \cite{MANIVEL}. For connections with the representation theory of the symmetric and general linear groups see Fulton \cite{BLUE} and Sagan \cite{SAGAN}. For more about $\lambda$-rings see Knutson \cite{KNUTSON}. For more about the $q$-exponential functions and $q$-binomial coefficients see the book by Kac and Cheung \cite{KAC} or the online notes by Foata and Han \cite{FOAT}. The standard reference for basic hypergeometric series is Gasper and Rahmen \cite{HYPGEOM}. For an interesting approach to obtaining basic hypergeometric series identities by specialization of symmetric functions see Bowman \cite{BOWMAN}

\subsubsection{Partitions}\label{partitions}

For our purposes, a composition is just a list of non-negative integers, finitely many of which are greater than zero. For most purposes we will identify lists of different lengths which differ only by a tail of zeros. 

A {\it partition} is a composition whose parts are weakly decreasing. If the sum of the parts of $\lambda$ is $n$ then we say that $\lambda$ is a partition of $n$ and write $\lambda \vdash n$. 

Some authors do not allow partitions to contain zero parts, however there are some circumstances in which it is convenient to envision the partition as sitting inside some rectangle of pre-determined dimensions, in which case the tail of zeros becomes significant.

\par\bigskip
The set of all partitions of $n$ with $k$ parts (any number of which may be zero) is denoted by $\mathcal{P}_k^n$. The generating series for partitions is the first $q$-exponential function:
$$ E_q(z) = \frac{1}{(z;q)_\infty} = \sum_{n,k} |\mathcal{P}_k^n| q^n z^k = \sum_k \frac{z^k}{(q;q)_k}$$
Here we are making use of the hypergeometric notation:
$$ (a;q)_n = \prod_{k=0}^{n-1} (1-aq^k) $$

It is common practice to represent a partition pictorially as rows of boxes. The {\it conjugate} of a partition is obtained by reflecting in the main diagonal. The prime symbol is used to indicate conjugation. Conjugating a partition interchanges the role of rows and columns.  The parts of $\lambda'$ correspond to the columns of $\lambda$. 

\par\bigskip

Below, on the left, is the diagram for the partition $\lambda = (4,2,1)$ and on the right, its conjugate $\lambda' = (3,2,1,1)$
$$ \tableau{{\ }&{\ }&{\ }&{\ }\\{\ }&{\ }\\{\ }} 
\,\,\,\,\,\,\,\,\,\,\,\,\,\,\,\,
\tableau{{\ }&{\ }&{\ }\\{\ }&{\ }\\{\ } \\{\ }}
$$ 
We are using the English convention.

\par\bigskip

We define the length of a partition, denoted by $l(\lambda)$, to be the number of nonzero parts. With this convention, the length of a partition is always less than or equal to its number of parts. For example, $\lambda = (3,2,2,1,0,0)$ is a partition of $8$ of length $4$ with $6$ parts. 

There are two natural additions on the set of partitions. The first one corresponds to the concatenation of the individual parts:
$$ \lambda + \mu = (\lambda_1 + \mu_1, \lambda_2 + \mu_2, \ldots, \lambda_n + \mu_n) $$
The second corresponds to taking the union of the two sets of parts, and then re-ordering them as appropriate. For example:
$$ (5,3,1) \cup (7,3,2,2) = (7,5,3,3,2,2,1) $$
One may check that:
$$ \lambda' + \mu' = (\lambda \cup \mu)' $$ 

\par\bigskip
Of particular interest are partitions that contain no repeated parts, including no repeated zero parts.
The smallest such partition with $k$ distinct parts is the staircase partition:
$$\delta_k = (k-1, k-2, \ldots, 2,1,0)$$ 

Observe that $\delta_k$ is a partition of $k(k-1)/2$. If $\lambda$ is a partition of $n$ into $k$ parts then $\lambda + \delta_k$ is a partition of $n + k(k-1)/2$ into $k$ distinct parts. Conversely, every partition of $n + k(k-1)/2$ into $k$ distinct parts is equal to $\lambda + \delta_k$ for some $\lambda \vdash n$.

\par\bigskip

We denote by $\mathcal{D}_k^n$ the set of all partitions of $n + k(k-1)/2$ into $k$ distinct parts (one of which may be zero).
The generating series for partitions with distinct parts is the second $q$-exponential function:
$$ e_q(z) = (-z;q)_\infty  = \sum_{n,k} |\mathcal{D}_k^n| q^n z^k = \sum_k q^{\binom{k}{2}} \frac{z^k}{(q;q)_k}  $$

Note that the conjugate of a partition with distinct parts will not, in general, be another partition with distinct parts. In fact, the only partition with distinct parts whose conjugate also has distinct parts is the staircase partition $\delta_k$.

\par\bigskip

It is perhaps worth pointing out here, though we won't need it until much later, that if $\pi_q$ is the function which maps $q$ to $1/q$ then:
$$ \pi_q[(q;q)_n] = (-1)^n q^{-\frac{n(n+1)}{2}} (q;q)_n$$
and so:
$$ \pi_q[E_q(z)] = e_q(-qz) $$

\par\bigskip
If the diagram for a partition $\mu$ sits properly inside the diagram for the partition $\lambda$, that is, if each part of $\mu$ is less than or equal to the corresponding part of $\lambda$ then we write $\mu \subseteq \lambda$. For such a pair of partitions, we may form the {\it skew partition} $\lambda / \mu$ which is simply the collection of boxes in $\lambda$ which are not also boxes of $\mu$.

A skew partiton $\lambda / \mu$ is said to be a {\it horizontal strip} if:
$$ \lambda_1 \geq \mu_1 \geq \lambda_2 \geq \ldots \geq \lambda_n \geq \mu_n $$
That is, no column contains more than a single box.
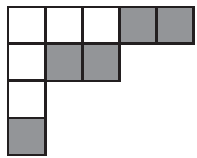
\[
\tableau{ &&&\graycell &\graycell \\ &\graycell &\graycell \\ \\ \graycell }
\]
For an arbitrary partition $\mu$ we let $U(\mu)$ denote the set of partitions $\lambda$ such that $\lambda / \mu$ is a horizontal strip, and let $D(\mu)$ denote the set of partitions $\nu$ such that $\mu / \nu$ is a horizontal strip. 

More specifically, we let $U_r(\mu)$ denote the set of partitions $\lambda$ such that $\lambda / \mu$ is a horizontal strip with exactly $r$ boxes, and let $D_r(\mu)$ denote the set of partitions $\nu$ such that $\mu / \nu$ is a horizontal strip with exactly $r$ boxes.

Similarly a skew partition $\lambda / \mu$ is a {\it vertical strip} if:
$$ \lambda'_1 \geq \mu'_1 \geq \lambda'_2 \geq \ldots \geq \lambda'_n \geq \mu'_n $$
That is, no {\it row} contains more than a single box.
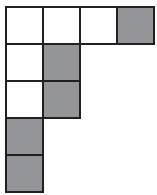
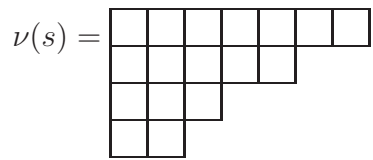
\[
\tableau{ &&&\graycell \\ &\graycell \\ &\graycell \\ \graycell \\ \graycell }
\]
For an arbitrary partition $\mu$ we let $\tilde{U}(\mu)$ denote the set of partitions $\lambda$ such that $\lambda / \mu$ is a vertical strip, and let $\tilde{D}(\mu)$ denote the set of partitions $\nu$ such that $\mu / \nu$ is a vertical strip. 

More specifically, we let $\tilde{U}_r(\mu)$ denote the set of partitions $\lambda$ such that $\lambda / \mu$ is a vertical strip with exactly $r$ boxes, and let $\tilde{D}_r(\mu)$ denote the set of partitions $\nu$ such that $\mu / \nu$ is a vertical strip with exactly $r$ boxes.

Clearly, if $\lambda / \mu$ is a horizontal strip, then the conjugate $\lambda' / \mu'$ is a vertical strip. 

\par\bigskip



An $(n,m)$-binomial path is a binary string containing exactly $n$ zeros and $m$ ones. Pictorially, we may represent an $(n,m)$-binomial path by drawing an $n$ by $m$ grid and tracing a path from the bottom left hand corner to the top right hand corner by reading the binary string from left to right and taking a step up each time we read a one,
and a step across each time we read a zero.

For example, if $s = 00110101000$ we have the following path:

\[
s=\vcenterbox{\mbox{\setlength{\unitlength}{\cellsize}
\begin{picture}(7,4)
\put(0,0){\line(1,0){7}}
\put(0,0){\line(0,1){4}}
\put(0,4){\line(1,0){7}}
\put(7,0){\line(0,1){4}}
\linethickness{2pt}
\put(0,0){\line(1,0){1}}
\put(1,0){\line(1,0){1}}
\put(2,0){\line(0,1){1}}
\put(2,1){\line(0,1){1}}
\put(2,2){\line(1,0){1}}
\put(3,2){\line(0,1){1}}
\put(3,3){\line(1,0){1}}
\put(4,3){\line(0,1){1}}
\put(4,4){\line(1,0){1}}
\put(5,4){\line(1,0){1}}
\put(6,4){\line(1,0){1}}
\end{picture}}}
\]

Observe that the boxes lying ``above'' a binomial path form a Young diagram. We shall call this diagram $\alpha(s)$.
Likewise, after a rotation, the boxes lying ``below'' a binomial path form a Young diagram. We shall call this diagram $\beta(s)$. In our example, we have:

\[
\alpha(s)=\tableau{\ &\ &\ &\ \\\ &\ &\ \\\ &\ \\\ &\ }
\,\,\,\,\,\,\,\,\,\,\,\,
\beta(s) = \tableau{\ &\ &\ &\ &\ \\\ &\ &\ &\ &\ \\\ &\ &\ &\ \\\ &\ &\ }
\]

If $s$ is an $(n,m)$ binomial path then the {\it dual} of $s$ is the $(n,m)$-binomial path
$\tilde{s}$ which would be obtained by reading the binary string of $s$ from right to left,
rather than from left to right. In our example $\tilde{s} = 00010101100$
\[
\tilde s=\vcenterbox{\mbox{\setlength{\unitlength}{\cellsize}
\begin{picture}(7,4)
\put(0,0){\line(1,0){7}}
\put(0,0){\line(0,1){4}}
\put(0,4){\line(1,0){7}}
\put(7,0){\line(0,1){4}}
\linethickness{2pt}
\put(0,0){\line(1,0){1}}
\put(1,0){\line(1,0){1}}
\put(2,0){\line(1,0){1}}
\put(3,0){\line(0,1){1}}
\put(3,1){\line(1,0){1}}
\put(4,1){\line(0,1){1}}
\put(4,2){\line(1,0){1}}
\put(5,2){\line(0,1){1}}
\put(5,3){\line(0,1){1}}
\put(5,4){\line(1,0){1}}
\put(6,4){\line(1,0){1}}
\end{picture}}}
\]

It is not hard to see that the upper diagram associated with a binomial path $s$
is equal to the lower diagram associated with the dual path $\tilde{s}$. 
Likewise the lower diagram of $s$ is equal to the upper diagram of $\tilde{s}$.
\[
\alpha(\tilde{s}) = \tableau{\ &\ &\ &\ &\ \\\ &\ &\ &\ &\ \\\ &\ &\ &\ \\\ &\ &\ } = \beta(s)
\, \, \,\,\,\,\,\,\,\,\,\,\,\,\,\,
\beta(\tilde{s}) = \tableau{\ &\ &\ &\ \\\ &\ &\ \\\ &\ \\\ &\ } = \alpha(s)
\]

Now let $s'$ denote the binary string obtained from $s$ by swapping all the zeros for ones and vice versa. In our case $s' = 11001010111$. Note that $s'$ is now an $(m,n)$ binomial path.
\[
s'=\vcenterbox{\mbox{\setlength{\unitlength}{\cellsize}
\begin{picture}(4,7)
\put(0,0){\line(1,0){4}}
\put(0,0){\line(0,1){7}}
\put(0,7){\line(1,0){4}}
\put(4,0){\line(0,1){7}}
\linethickness{2pt}
\put(0,0){\line(0,1){1}}
\put(0,1){\line(0,1){1}}
\put(0,2){\line(1,0){1}}
\put(1,2){\line(1,0){1}}
\put(2,2){\line(0,1){1}}
\put(2,3){\line(1,0){1}}
\put(3,3){\line(0,1){1}}
\put(3,4){\line(1,0){1}}
\put(4,4){\line(0,1){1}}
\put(4,5){\line(0,1){1}}
\put(4,6){\line(0,1){1}}
\end{picture}}}
\]
Observe that $(\alpha(s'), \beta(s')) = (\beta(s)', \alpha(s)')$, and similarly $(\alpha(\tilde{s}'), \beta(\tilde{s}')) = (\beta(\tilde{s})', \alpha(\tilde{s})') = (\alpha(s)', \beta(s)')$. In other words, reversing the string and interchanging the role of zeros and ones is equivalent to conjugating the partition.

\[
\alpha(s') =
\tableau{\ &\ &\ &\ \\\ &\ &\ &\ \\\ &\ &\ &\ \\\ &\ &\ \\\ &\ } = \beta(s)'
\, \, \,\,\,\,\,\,
\beta(s') =\tableau{\ &\ &\ &\ \\\ &\ &\ &\ \\\ &\ \\\ } = \alpha(s)'
\]

\par\bigskip

There is a also a way that we may associate a pair of diagrams $(\nu(s), \gamma(s))$ with {\it distinct} parts
to any given binomial path. Suppose that we label the steps of $s$ from left to right with the integers from $0$ to $n + m -1$.

\[
s=\vcenterbox{\mbox{\setlength{\unitlength}{\cellsize}
\begin{picture}(7,4)
\put(0,0){\line(1,0){7}}
\put(0,0){\line(0,1){4}}
\put(0,4){\line(1,0){7}}
\put(7,0){\line(0,1){4}}
\linethickness{2pt}
\put(0,0){\line(1,0){1}}
\put(1,0){\line(1,0){1}}
\put(2,0){\line(0,1){1}}
\put(2.1,0.2){2}
\put(2,1){\line(0,1){1}}
\put(2.1,1.2){3}
\put(2,2){\line(1,0){1}}
\put(3,2){\line(0,1){1}}
\put(3.1,2.2){5}
\put(3,3){\line(1,0){1}}
\put(4,3){\line(0,1){1}}
\put(4.1,3.2){7}
\put(4,4){\line(1,0){1}}
\put(5,4){\line(1,0){1}}
\put(6,4){\line(1,0){1}}
\end{picture}}}
\]

Then $\nu(s)$ is the diagram which has a row of length $k$ if and only if there is some upstep
of $s$ labelled with $k$. We have of course that $\nu(s) = \alpha(s) + \delta_n$ 

\[
\nu(s)=
\tableau{\ &\ &\ &\ &\ &\ &\ \\\ &\ &\ &\ &\ \\\ &\ &\ \\\ &\ }
\]

The other partition $\gamma(s)$ is the diagram which has a row of length $k$ if and only if there is some
{\it across} step which is labelled with a $k$. 
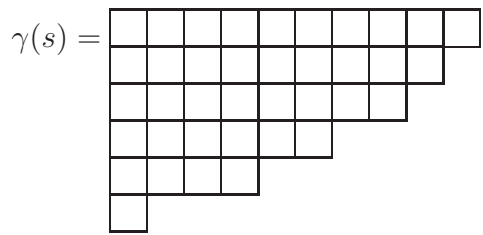
\[
\gamma(s)=
\tableau{\ &\ &\ &\ &\ &\ &\ &\ &\ &\ \\\ &\ &\ &\ &\ &\ &\ &\ &\ \\\ &\ &\ &\ &\ &\ &\ &\ \\\ &\ &\ &\ &\ &\ \\\ &\ &\ &\ \\\ }
\]
Of course $\nu(s) \cup \gamma(s) = \delta_{m+n}$. 

\par\bigskip
Performing the same procedure with the dual path $s'$ we get:

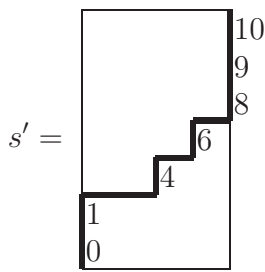
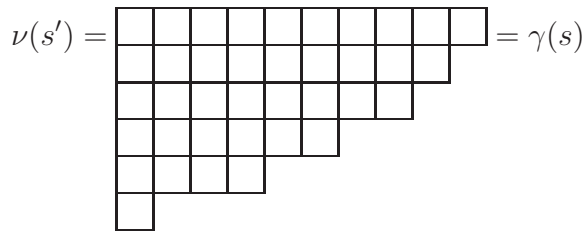
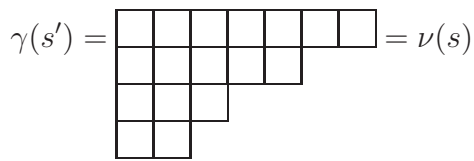
\[
s'=\vcenterbox{\mbox{\setlength{\unitlength}{\cellsize}
\begin{picture}(4,7)
\put(0,0){\line(1,0){4}}
\put(0,0){\line(0,1){7}}
\put(0,7){\line(1,0){4}}
\put(4,0){\line(0,1){7}}
\linethickness{2pt}
\put(0,0){\line(0,1){1}}
\put(0.1,0.2){0}
\put(0,1){\line(0,1){1}}
\put(0.1,1.2){1}
\put(0,2){\line(1,0){1}}
\put(1,2){\line(1,0){1}}
\put(2,2){\line(0,1){1}}
\put(2.1,2.2){4}
\put(2,3){\line(1,0){1}}
\put(3,3){\line(0,1){1}}
\put(3.1,3.2){6}
\put(3,4){\line(1,0){1}}
\put(4,4){\line(0,1){1}}
\put(4.1,4.2){8}
\put(4,5){\line(0,1){1}}
\put(4.1,5.2){9}
\put(4,6){\line(0,1){1}}
\put(4.1,6.2){10}
\end{picture}}}
\]
This time:
\[
\nu(s')=
\tableau{\ &\ &\ &\ &\ &\ &\ &\ &\ &\ \\\ &\ &\ &\ &\ &\ &\ &\ &\ \\\ &\ &\ &\ &\ &\ &\ &\ \\\ &\ &\ &\ &\ &\ \\\ &\ &\ &\ \\\ }  = \gamma(s)
\]
and:
\[
\gamma(s')=
\tableau{\ &\ &\ &\ &\ &\ &\ \\\ &\ &\ &\ &\ \\\ &\ &\ \\\ &\ } = \nu(s)
\]

Observe that, a little surprisingly, we have $\gamma(s) = \nu(s')$ and $\nu({s}) = \gamma(s')$, which gives us the following proposition:
\begin{proposition} \label{flip} Let $\lambda$ be any partition contained within an $n$ by $m$ box, and let $\mu$ be its compliment with respect to this box, then we have:
$$ (\lambda + \delta_n) \cup (\mu' + \delta_m) = \delta_{n + m} $$
\begin{proof}
Suppose that $\lambda = \alpha(s)$, so that $\lambda + \delta_n = \nu(s)$. Now $\mu' = \alpha(s')$ and $\mu' + \delta_m = \nu(s') = \gamma(s)$. The result now follows from the fact that $\nu(s)\cup \gamma(s) = \delta_{n+m}$.
\end{proof}

\end{proposition}

\subsubsection{Monomial symmetric functions}
Consider now the multivariate polynomial ring $\mathbb{Q}[x_1, \ldots, x_k]$. 
Each monomial in $\mathbb{Q}[x_1, \ldots, x_k]$ corresponds to a composition. For example, the monomial $x_1 x_3^2$ in $\mathbb{Q}[x_1,x_2,x_3]$ corresponds to the list $(1,0,2)$.

For notational convenience,
if $\eta = (\eta_1, \eta_2, \ldots, \eta_k)$ is a composition, then by $X^\eta$ we mean the monomial $x_1^{\eta_1} x_2^{\eta_2} \ldots x_k^{\eta_k}$. There is a natural addition on the space of compositions which corresponds to multiplication in the polynomial ring. if $\eta = (\eta_1, \ldots, \eta_k)$ and $\gamma = (\gamma_1,  \ldots, \gamma_k)$ then $\eta + \gamma = (\eta_1 + \gamma_1,  \ldots, \eta_k + \gamma_k)$, and $X^{(\eta + \gamma)} = X^\eta X^\gamma$.

\par\bigskip

The symmetric group $\mathcal{S}_k$ acts on the set of compositions with $k$ parts  by permuting the parts. Each orbit of $\mathcal{S}_k$ contains a unique partition. 
For any composition $\lambda$ let $r(\lambda,i)$ denote the number of parts of $\lambda$ equal to $i$. Let us define:
$$r_\lambda = \prod_{i \geq 0} r(\lambda,i)!$$
Then $r_\lambda$ is the order of the subgroup of $\mathcal{S}_k$ that stabilizes the composition ${\lambda}$. If $\lambda$ is a partition with distinct parts then $r_\lambda = 1$.
 
\par\bigskip

For each partition $\lambda$ we define the {\it monomial symmetric polynomial} to be:
$$ m_\lambda = \frac{1}{r_\lambda} \sum_{\sigma \in S_n} X^{\sigma(\lambda)} $$

Since $r_\lambda$ is the number of permutations in $\mathcal{S}_k$ that will stabilize the monomial
$x^\lambda$, the monomial symmetric
function is the sum of all {\it distinct} permutations of this monomial.


As $\lambda$ runs over the set of all permutations of $n$ with at most $k$ parts, the set $\{m_{\lambda} | \lambda \in \mathcal{P}^k_n \}$ forms a basis for the vector space $\Lambda_k^n$ consisting of all symmetric polynomials of homogeneous degree $k$ in $n$ variables. We shall write:
$$ \Lambda_n = \bigoplus_{k} \Lambda_n^k$$
to denote the ring of symmetric functions in $n$ variables. 
The map:
$$ \pi_n : \Lambda_{n} \to \Lambda_{n-1}$$
which sets the $n$th variable equal to zero is a homomorphism whose kernel is generated by the monomial symmetric functions indexed by partitions with exactly $n$ nonzero parts. 

\subsubsection{Plethystic notation}\label{pleth}
It is often more convenient to work in infinitely many variables. It is also convenient to make use of the plethystic notation. In the plethystic notation one writes:
$ X = x_1 + x_2 + x_3 + \cdots $
to denote the set of indeterminates:
$ X = \{ x_1, x_2, x_3, \ldots \}$
and $ Y = y_1 + y_2 + y_3 + \cdots$
to denote the set of indeterminates
$ Y = \{ y_1, y_2, y_3, \ldots \}$

Furthermore by $X+Y$ we denote the disjoint union of the variables in the sets $X$ and $Y$ 
\begin{align*}
X + Y & =  x_1 + x_2 + \cdots + y_1 + y_2 + \cdots \\
& =  \{x_1, x_2, \ldots, y_1, y_2, \ldots \} \\
& =  X \cup Y
\end{align*}
and by $XY$ we denote the cartesian product:
\begin{align*}
XY & =  (x_1 + x_2 + \cdots)(y_1 + y_2 + \cdots) \\
   & =  x_1 y_1 + x_1 y_2  + \cdots + x_2 y_1 + x_2 y_2 + \cdots \\
   & =  \{x_1 y_1, x_1 y_2, \ldots x_2 y_1, x_2 y_2, \ldots \} \\
   & =  X \times Y
\end{align*}
The complete and elementary symmetric functions may be defined in terms of their generating functions:
$$ \Omega_z(X) = \Omega(Xz) = \prod_{x \in X} \frac{1}{1-xz} = \sum_{n=0}^\infty h_n(X) z^n $$
$$ \tilde{\Omega}_z(X) = \tilde{\Omega}(Xz) = \prod_{x \in X}(1+xz) = \sum_{n=0}^\infty e_n(X) z^n $$
Note that the first $q$-exponential funtion may be expressed as:
$$ E_q(z) = \Omega\left ( \frac{z}{1-q} \right ) $$
where $\frac{1}{1-q}$ is the alphabet $1 + q + q^2 + \dots $.
Similarly the second $q$-exponential function may be expressed as:
$$ e_q(z) = \tilde{\Omega}\left ( \frac{z}{1-q} \right ) $$
\noindent
The involution $\omega$ is defined by:
$$ \omega[h_n(X)] = e_n(X) $$

(extended multiplicatively). The ``forgotten'' symmetric functions are defined by:
$$ \omega[m_\lambda(X)] = f_\lambda(X) $$

\par\bigskip
For each variable $x$ let us define an {\it anti-variable} $\overline{x}$ such that:
$$ \{x\} \cup \{ \overline{x} \} = \{ \}$$
Alternatively one can think in terms of {\it anti-sets}:
$$ \{x\} \cup \overline{ \{ x \} } = \{ \}$$
Now let
$-X = \overline{x}_1 + \overline{x}_2 + \overline{x}_3 + \cdots $
denote the set of anti-variables
$ \{ \overline{x}_1, \overline{x}_2,\overline{x}_3, \ldots \} $
or, the {\it anti-set} of variables
$ \overline{\{x_1, x_2, x_3, \ldots \}}$.

\par\bigskip
For an arbitrary symmetric function $f$ we define:
$$ f(-X) = \omega f(\epsilon X) $$
where $\epsilon X$ is the more obvious but less natural form of negation:
\begin{align*}
\epsilon X & = \{-x_1, -x_2, -x_3, \ldots \}
\end{align*}
We may now express the elementary symmetric functions in terms of the complete symmetric functions:
$$ \tilde{\Omega}(X) = \Omega(-\epsilon X) $$

\subsubsection{Schur functions}\label{schurfun}
The most important basis for the ring of symmetric functions, from the perspective of geometry or representation theory,  is the Schur functions.
 
Let $H(X)$ denote the infinite Toepliz matrix $(h_{i-j}(X))$, where we are assuming that $h_{-k}(X) = 0 $ for all $k \geq 1$. In other words the matrix is upper triangular with ones on the diagonal. For $\lambda$ and $\mu$ partitions with $n$ parts (any number of which may be zero), the skew Schur function $s_{\lambda / \mu}(X)$ may be defined by minors of this matrix as follows:
$$ s_{\lambda / \mu}(X) = \det(H_{I,J}(X)) = \det(h_{\lambda_i - \mu_j + j - i}(X))$$
where $I = \mu + \delta_{n}$
and $J = \lambda + \delta_{n}$. Note that we indexing both the rows and the columns from zero rather than one. In particular the regular Schur functions are given by:
$$ s_\lambda(X) = \det(H_{[n] \, , \, \lambda + \delta_n}) = \det(h_{\lambda_i + j - i}(X)) $$
where $[n] = \{0,1,2,...,n-1\}$. Note that $n$ does not appear explicitly in this equation, and if $n$ is taken to be larger than the number of nonzero parts of $\lambda$ then the resulting matrix is block diagonal, with the first block of size $l(\lambda)$ independent of $n$, and the remaining blocks containing ones, which do not affect the determinant.

The fact that the generating function for the complete symmetric functions has the property
$ \Omega_z(X + X') = \Omega_z(X) \Omega_z(X') $
lifts to the fact that:
$ H(X + X') = H(X)H(X') $,
and so, by the Cauchy--Binet theorem we have:
$$ s_\lambda(X + X') = \sum_{\mu \subseteq \lambda} s_{\lambda / \mu}(X) s_\mu(X') $$
Consider now the infinite matrix $E(X) = (e_{i-j}(X))$
The fact that $ \Omega(X)\tilde{\Omega}(\epsilon X)  = 1 $ now lifts to the fact that $H(X)E( \epsilon X) = I$. Thus by Jacobi's formula for the minors of the inverse matrix, we have:
\begin{align*}
\det(H_{I,J}(X)) & =  (-1)^{\sum_{k=1}^r j_k - i_k} \det(E_{J',I'}(\epsilon X)) \\
& =  \det(E_{J',I'}(X))
\end{align*}
Now by proposition \ref{flip}, if these are $(n+m)$ by $(n+m)$ matrices and $I = \lambda + \delta_n$ then $I' = \lambda' + \delta_m$. In other words we have an alternative expression for the Schur functions given by:
\begin{align*}
s_{\lambda / \mu}(X) 
& =  \det( e_{\lambda'_i - \mu'_j + j - i}(X))
\end{align*}
We also get from this that:
$$ \omega[s_\lambda(X)] = s_{\lambda'}(X) $$

\subsection{The Umbral Calculus}\label{umbralcalc}

A nice reference for finite dimensional linear algebra is Hoffman and Kunze \cite{LIN}. For more about algebras, co-algebras, Hopf algebras and quantum groups see the online notes by Arun Ram \cite{RAM}. The book by Roman \cite{ROMAN} discusses the infinite-dimensional subtleties more carefully than is done here.

\subsubsection{Coalgebras}

\noindent
An algebra is a vector space $V$ over some field $F$ equipped with a multiplication:
$$ m : V \otimes V \to V$$
and a unit:
$$ e : F \to V $$
such that the following two diagrams commute:
\begin{diagram}
V\otimes V\otimes V &\rTo^{m\otimes id} & V\otimes V\\
\dTo^{id\otimes m}&& \dTo_{m}\\
V\otimes V &\rTo^{m}&V
\end{diagram}

\par\bigskip
This is just the associative law.

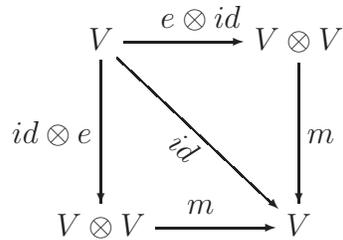
\begin{diagram}
V &\rTo^{e\otimes id} & V\otimes V\\
\dTo^{id\otimes e}& \rdTo_{id} & \dTo_{m}\\
V\otimes V &\rTo^{m}&V
\end{diagram}


\noindent
This says that multiplication by the identity has no effect.

\par\bigskip
A commutative algebra is an algebra with the additional property that:
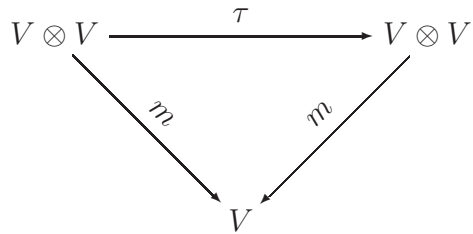
\begin{diagram}
V\otimes V &&\rTo^{\tau} && V\otimes V\\
&\rdTo^{m} &&\ldTo^{m}&\\
&&V&&
\end{diagram}
\noindent
where $\tau : V \otimes V \to V \otimes V$ is the map
$ \tau( x \otimes y ) = y \otimes x $.

The ring $\Lambda$ of symmetric functions, and the ring  $\mathbb{Q}[x]$ of polynomials in one variable are both commutative algebras.

An algebra morphism is a vector space morphism:
$ \psi : V \to V $
satisfying the two properties:
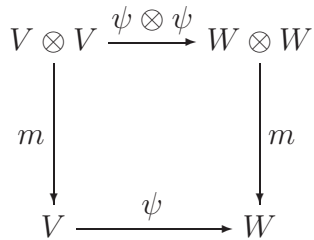
\begin{diagram}
V\otimes V &\rTo^{\psi\otimes\psi} & W\otimes W\\
\dTo^{m}&& \dTo_{m}\\
V &\rTo^{\psi}&W
\end{diagram}
The multiplication is preserved.

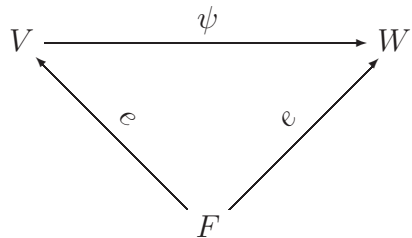
\begin{diagram}
V &&\rTo^{\psi} && W\\
&\luTo^{e} &&\ruTo^{e}&\\
&&F&&
\end{diagram}
The identity is preserved.

\par\bigskip
The map which sends the complete symmetric function $h_n(X)$ to $x^n$ is an algebra morphism from $\Lambda$ onto $\mathbb{Q}[x]$.

\par\bigskip
A coalgebra is a vector space $V$ equipped with a comultiplication:
$$ \Delta : V \to V \otimes V $$
and a counit:
$$ c : V \to F $$
such that the following diagrams commute:
\begin{diagram}
V\otimes V\otimes V &\lTo^{\Delta\otimes id} & V\otimes V\\
\uTo^{id\otimes \Delta}&& \uTo_{\Delta}\\
V\otimes V &\lTo^{\Delta}&V
\end{diagram}
This is the co-associative law.
\begin{diagram}
V &\lTo^{c\otimes id} & V\otimes V\\
\uTo^{id\otimes c}& \ruTo_{id} & \uTo_{\Delta}\\
V\otimes V &\lTo^{\Delta}&V
\end{diagram}
This says that co-multiplying by the co-unit has no affect.

\par\bigskip
A co-commutative co-algebra is a co-algebra with the additional property that:
\begin{diagram}
V\otimes V &&\lTo^{\tau} && V\otimes V\\
&\luTo^{\Delta} &&\ruTo^{\Delta}&\\
&&V&&
\end{diagram}
\noindent
Both the ring of symmetric functions and the ring of polynomials are co-commutative co-algebras. In the case of $\mathbb{Q}[x]$ we have:
$$ \mathbb{Q}[x] \otimes \mathbb{Q}[x] \cong \mathbb{Q}[x,x']$$
and the coproduct is given by:
$$\Delta[p(x)] = p(x + x')$$
In the symmetric function case we have, when working in an infinite alphabet, that:
$$ \Lambda \otimes \Lambda \cong \Lambda $$
The coproduct is given by the plethystic addition of alphabets.

\par\bigskip
A co-algebra morphism is a vector space morphism:
$ \phi : V \to V $
satisfying the two properties:
\begin{diagram}
V\otimes V &\lTo^{\phi\otimes\phi} & W\otimes W\\
\uTo^{\Delta}&& \uTo_{\Delta}\\
V &\lTo^{\phi}&W
\end{diagram}
The comultiplication is preserved.
\begin{diagram}
V &&\lTo^{\phi} && W\\
&\rdTo^{c} &&\ldTo^{c}&\\
&&F&&
\end{diagram}
The counit is preserved.

The map which sends $x^n/n!$ to the complete symmetric function $h_n(X)$ is a co-algebra morphism from $\mathbb{Q}[x]$ into $\Lambda$. We shall make use of this homomorphism to lift properties of $\mathbb{Q}[x]$ to $\Lambda$.

\par\bigskip

Categorically, the notions of algebra and co-algebra are dual. Combinatorially the product describes how the elements of an algebra can be combined, while the coproduct describes how the elements of a co-algebra may be decomposed \cite{ROTA-COALG1, ROTA-COALG2}.

\par\bigskip
Every vector space $V$ has a dual $V^*$ which is the space of linear functionals from $V$ into the field $F$. There is a natural map:
$$ ( - , - ): V^* \otimes V \to F $$
given by:
$$ ( w, v ) = w[v] $$

\par\bigskip
In finitely many dimensions, every vector space is isomorphic to its dual, and once an isomorphism $\phi : V \to V^*$ has been fixed we may define a non-degenerate bilinear form:
$$ \langle - , - \rangle : V \otimes V \to F $$
by:
$$ \langle v_1, v_2 \rangle = (\phi(v_1), v_2 ) $$

\par\bigskip
Still in finite dimensions, if $V$ is an algebra then $V^*$ is a coalgebra with coproduct:
$$ ( \Delta w, v_1 \otimes v_2 ) = ( w, m(v_1, v_2) ) $$
and counit:
$$c(w) = (w, e(1))$$

\par\bigskip
Furthermore, if $\psi : V \to W$ is an algebra morphism, then its adjoint $\psi^* : W^* \to V^*$ is a coalgebra morphism.

\bigskip
\noindent
Similarly, if $V$ is a co-algebra then $V^*$ is an algebra with multiplication: 
$$ ( m(w_1, w_2), v ) = ( w_1 \otimes w_2, \Delta(v) ) $$
and unit:
$$ (e(x), v) = x c(v)  $$
If $\phi : V \to V$ is an co-algebra morphism, then its adjoint $\phi^* : V^* \to V^*$ is an algebra morphism.

\par\bigskip
Unfortunately, both $\Lambda$ and $\mathbb{Q}[x]$ are infinite dimensional. In infinitely many dimensions, a vector space $V$ need not be isomorphic to its dual $V^*$, and instead of equality we have:
$$ V^* \otimes V^* \subseteq (V \otimes V)^* $$

Since $\Delta(w) \in V^* \otimes V^*$ whenever $w \in V^*$ we are still able, in infinite dimensions, to induce a product structure on the dual from an existing product structure on the original space, and we still have that algebra morphisms of $V^*$ correspond to co-algebra morphisms of $V$. The converse is however no longer true.





\subsubsection{Sequences of Binomial Type}\label{bintype}
The classification of the co-algebra isomorphisms of the polynomial ring $\mathbb{Q}[x]$ by consideration of the algebra morphisms of the dual space goes back to Gian-Carlo Rota \cite{ROTA,FINOP}. Since then there have been many other nice expositions such as Gessel \cite{GESSUMB} and Roman \cite{ROMAN}, as well as attempts to find $q$-analogs \cite{ANDREWS, QUMBRAL, QROMAN}, and attempts to generalize the umbral calculus to symmetric functions: Loeb \cite{LOEB}, Mendez \cite{COMBSYM1, MENDEZ} and Chen \cite{CHEN2, CHEN1}.  We follow here, in particular, the presentation of Garsia \cite{EXPOSE, UMBLAG, COMPSEQ}. See the following online resources \cite{UMBSURV1, UMBSURV2} for a more extensive bibliography.

The ring of polynomials in one indeterminate $V = \mathbb{Q}[x]$ has a distinguished Hamel basis given by:
$$\{1,x,\frac{x^2}{2!},\frac{x^3}{3!}, \ldots \}$$ 
The dual space $V^* = \mathbb{Q}[[y]]$ has a distinguished Schauder basis:
$$\{1,y,y^2,y^3, \ldots \}$$
We may think of $y^k$ as the $kth$ coefficient extraction operator (acting on exponential generating series):
$$ (y^k, p(x) ) = \left [ \frac{x^k}{k!} \right ] p(x) $$
Note that coefficient extraction and differentiation are related by:
$$ \left [ \frac{x^k}{k!} \right ]f(x) = \mathcal{LD}^k f(x) $$
where $\mathcal{D}$ denotes the ordinary differential operator $d/dx$, and $\mathcal{L}$ is the constant term operator. Thus we may also write the pairing between $V$ and $V^*$ as:
$$ (f(y), p(x) ) = \mathcal{L}f(\mathcal{D})[p(x)] $$

\noindent
By the Leibniz rule, the co-multiplication on $V$:
$$ \Delta p(x) = p(x + x') $$
induces the multiplication on $V^*$:
$$m(y^r,y^s) = y^{r+s}$$
Define a delta series to be a formal power series $f(y)$ which has a compositional inverse. For every delta series $f(y)$ the map:
$$ \phi_f(y) = f(y) $$ 
induces an algebra morphism of $V^*$ with inverse:
$$ \phi_g(y) = g(y) $$
where $g(y)$ is the compositional inverse of $f(y)$, that is: 
$$f(g(y)) = g(f(y)) = y$$ 

\noindent
A sequence of polynomials $ \{ p_n(x)\}$ is said to be of {\it binomial type} if:
$$ p_n(x + x') = \sum_k \binom{n}{k} p_k(x) p_{n-k}(x') $$
A sequence of polynomials $\{ p_n(x)\}$ is said to be of {\it convolution type} if:
$$ p_n(x + x') = \sum_k p_k(x) p_{n-k}(x') $$
Clearly, if $p_n(x)$ is of binomial type, then $p_n(x) / n!$ is of convolution type. 

\par\bigskip
Sequences of convolution type arise as images, under co-algebra maps, of the standard basis $\left \{ x^n / n! \right \}$ while the corresponding sequences of binomial type arise as images, under co-algebra maps, of the basis $\{x^n\}$.

By duality, the co-algebra maps of $V$ are each adjoint to an algebra map of $V^*$. 
We shall write $ U_f = \phi_f^* $ to denote the umbral operator defined by:
$$U_f \left [\frac{x^n}{n!} \right ] = p_n(x)$$ 
where $\{ p_n(x)\}$ is the convolution type basis of $V$ which is dual to the basis $ \{1, g(y), g(y)^2, \ldots \} $ of $V^*$. As before $g(y)$ is the compositional inverse of $f(y)$.

\par\bigskip
Let $f_1(y) = \sum_k a_k y^k$ and $f_2(y) = \sum_k b_k y^k$ be delta series. Now let:
\begin{align*}
\alpha_{nk} & =  [y^n] f_2(y)^k \\
\beta_{nk} & =  [y^n] f_1(y)^k \\
\gamma_{nk} & =  [y^n] f_1(f_2(y))^k 
\end{align*}
and let $A = (\alpha_{nk})$ and $B = (\beta_{nk})$ and $C = (\gamma_{nk})$ be infinite matrices with $n,k \geq 1$.
Since:
\begin{align*}
f_1(f_2(y))^n & =  \sum_i \beta_{in} f_2(y)^i \\
& =  \sum_k \left (\sum_i \alpha_{ki} \beta_{in} \right)  y^k
\end{align*}
We have that:
$$AB = C$$
These are called Jabotinsky matrices \cite{COMPSEQ}.

\begin{proposition}
The umbral operator has the explicit expression:
$$ U_f = \sum_n \frac{x^n}{n!} \mathcal{L}f^n(\mathcal{D}) $$
Alternatively, the umbral operator may be characterized by the two properties:
$$ \mathcal{D} \circ U_f = U_f \circ g(\mathcal{D}) \leqno{(i)}$$
$$ \mathcal{L} \circ U_f = \mathcal{L} \leqno{(ii)} $$
\begin{proof}
For the first part, note that:
$$ [y^n] f(y) = f(\mathcal{D})\left[ \frac{x^n}{n!} \right ] $$
Let $A = \{a_{nk}\}$ be the matrix:
$$ a_{nk} = [y^k] f(y)^n $$
and let $B = \{b_{nk}\}$ be the inverse matrix:
$$ b_{nk} = [y^k] g(y)^n $$
We have:
\begin{align*}
\langle g(y)^k, U_{f}[x^n/ n!] \rangle
& =  \mathcal{L}g(\mathcal{D})^k \sum_j \frac{x^j}{j!} [y^n] f(y)^j \\
& =  \sum_j [y^j] g(y)^k [y^n] f(y)^j \\
& =  \sum_j b_{kj} a_{jn} \\
& =  \delta_{nk}
\end{align*}
For the second part, we have:
\begin{align*}
\langle g(y)^k, g(\mathcal{D})[p_n(x)] \rangle
& =  \mathcal{L}g(\mathcal{D})^{k+1}[p_n(x)] \\
& =  \delta_{k+1,n} \\
& =  \delta_{k, n-1} \\
& =  \langle g(y)^k, p_{n-1}(x) \rangle 
\end{align*}
Which shows that the umbral operator satisfies the first property. The second property is obvious. To see that these properties characterize the umbral operator, recall Taylor's formula:
$$ \mathcal{I} = \sum_n \frac{x^n}{n!} \mathcal{LD}^n $$
Suppose that $\mathcal{U}$ is any operator satisfying the above two properties, then we have:
\begin{align*}
\mathcal{U} & =  \sum_n \frac{x^n}{n!} \mathcal{LD}^n \mathcal{U} \\
& =  \sum_n \frac{x^n}{n!} \mathcal{L} \mathcal{U} \circ f^n(\mathcal{D}) \\
& =  \sum_n \frac{x^n}{n!} \mathcal{L} \circ f^n(\mathcal{D}) \\
& =  \mathcal{U}_f
\end{align*}
\qedhere
\end{proof}
\end{proposition}

\subsection{The Hall inner-product}

\par\bigskip

The results in this section are taken from the book of Lascoux \cite{LASCOUX} and the two papers by Garsia, Haiman and Tesler \cite{GARSIA2, GARSIA1}. See also the paper by Zabrocki \cite{ZABROCKI}.

\subsubsection{Preliminaries}
The Hall inner product may be defined by:
$$ \langle m_\lambda(Y), h_\mu(X) \rangle = \delta_{\lambda \mu} $$
Normally this is thought of as a map:
$$ \langle - , - \rangle : \Lambda \otimes \Lambda \to \mathbb{Q} $$
but we want to think of it as a map:
$$ \langle - , - \rangle : \Lambda^* \otimes \Lambda \to \mathbb{Q} $$

In analogy with the one variable case, $\Lambda$ is symmetric polynomials while $\Lambda^*$ is symmetric ``formal power series''. In lieu of the fact that a pair of bases $\{P_\lambda(X)\}$ and $\{Q_\lambda(Y)\}$ are dual with respect to the Hall inner product if and only if:
$$ \sum_\lambda P_\lambda(X) Q_\lambda(Y) = \Omega(XY) = \prod_{\substack{x \in X \\ x \in Y}} \frac{1}{1-xy}	$$
We shall always write elements of $\Lambda$ as symmetric functions in the alphabet $X$ and elements of $\Lambda^*$ as symmetric functions in the alphabet $Y$.

\par\bigskip

Any symmetric function $f(X)$ may be thought of as the operator ``multiplication by $f(X)$'', and as such there is an adjoint operator, which we denote by $\partial_{f}$. By definition:
$$ \langle \partial_{f}[g(Y)], h(X) \rangle = \langle g(Y), f(X) h(X) \rangle $$

The next proposition says that, via the Hall inner product, the coproduct on $\Lambda$ (given by plethystic addition of alphabets) induces the natural product on the dual. In fact, this property characterizes the Hall inner product.

\begin{proposition} \label{dual}
If $\{P_\lambda(X)\}$ and $\{Q_\lambda(Y)\}$ are a pair of dual bases (with respect to the Hall inner product), and:
$$ P_\lambda(X + X') = \sum_{\mu, \nu} d^{\lambda}_{\mu \nu} P_\nu(X) P_\mu(X') $$
then:
$$ Q_\mu(Y) Q_\nu(Y) = \sum_\lambda d^{\lambda}_{\mu \nu} Q_\lambda(Y)$$

\begin{proof}
Let:
$$ P_\lambda(X + X') = \sum_{\mu, \nu} d^{\lambda}_{\mu \nu} P_\nu(X) P_\mu(X') $$
and let:
$$ Q_\mu(Y) Q_\nu(Y) = \sum_\lambda \tilde{d}^{\lambda}_{\mu \nu} Q_\lambda(Y)$$
We have:
\begin{align*}
\Omega((X + X')Y) & =  \sum_\lambda P_\lambda(X + X') Q_\lambda(Y) \\
& =  \sum_\lambda Q_\lambda(Y) 
\sum_{\gamma, \mu} d^{\lambda}_{\gamma,\mu} P_\gamma(X) P_\mu(X') 
\\
\end{align*}
while:
\begin{align*}
\Omega(XY)\Omega(X'Y) & =  \left ( \sum_\gamma P_\gamma(X) Q_\gamma(Y) \right ) \left ( \sum_\mu P_\mu(X') Q_\gamma(Y) \right ) \\
& =  \sum_{\gamma, \mu} P_\gamma(X) P_\mu(X') Q_\gamma(Y) Q_\mu(Y) \\
& =  \sum_{\gamma, \mu} P_\gamma(X) P_\mu(X') \sum_\lambda \tilde{d}^{\lambda}_{\mu \gamma} Q_\lambda(Y) \\
& =  \sum_\lambda Q_\lambda(Y) \sum_{\gamma, \mu} \tilde{d}^{\lambda}_{\gamma \mu} P_\gamma(X) P_\mu(X')
\end{align*}
The result follows now since $\Omega(XY + X'Y) = \Omega(XY)\Omega(X'Y)$.
\end{proof}
\end{proposition}

\begin{corollary}[Taylor's theorem]
If $\{P_\lambda(X)\}$ and $\{ Q_\lambda(Y) \}$ are a pair of dual bases then for any symmetric function $f(X)$ we have:
$$ f(X + X') = \sum_\gamma \partial_{Q_\gamma}[f(X)] P_\gamma(X') $$
\begin{proof}
By linearity, it is sufficient to check on the basis $\{P_\lambda(X)\}$, from which the proposition follows by observation that:
\begin{align*}
\langle Q_\mu(Y), \partial_{Q_\gamma}[P_\lambda(X)] \rangle
& =  \langle Q_\mu(Y) Q_\gamma(Y), P_\lambda(X) \rangle \\
& =  d^{\lambda}_{\mu \gamma}
\end{align*}
\qedhere
\end{proof}
\end{corollary}
\noindent
Note that as a consequence of Taylor's theorem, and the self-duality of Schur functions, the fact that:
$$ s_\lambda(X + X') = \sum_{\mu} s_{\lambda / \mu}(X) s_\mu(X')$$
implies:
$$ s_{\lambda / \mu}(X) = \partial_{s_\mu}[s_\lambda(X)] $$

\subsubsection{Column operators}
Let us introduce now the ``translation'' and ``multiplication'' operators of Garsia, Haiman and Tesler \cite{GARSIA2, GARSIA1}:
\begin{equation} \label{trans}
\mathcal{T}(z)[ f(X) ] = f \left( X + z \right ) 
\end{equation}
\begin{equation} \label{mult}
\mathcal{P}(z)[ f(X) ] = \Omega(Xz)f(X) 
\end{equation}

\par\bigskip

Working with the pair of dual bases $\{m_\lambda(Y)\}$ and $\{h_\lambda(X)\}$, and noting that unless $\lambda$ is a row partition, $m_\lambda(X)$ will vanish when $X$ is an alphabet containing a single letter, we get, by Taylor's theorem, that:
$$ \mathcal{T}(z) = \sum_n z^n \partial_{h_n} $$
Similarly:
$$ \mathcal{P}(z) = \sum_n z^n h_n(X) $$
Clearly we have:
\begin{equation}
\mathcal{T}(z) = \mathcal{P}^*(z) 
\end{equation}

\par\bigskip

Similarly, working with the dual bases $\{f_\lambda(Y) \}$ and $\{e_\lambda(X)\}$ and recalling that $f_\lambda(-Y) = m_\lambda(\epsilon Y) $ we get that:
$$ \mathcal{T}(-z) = \sum_n (\epsilon z)^n \partial_{e_n} $$
and:
$$ \mathcal{P}(-z) = \sum_n (\epsilon z)^n e_n(X) $$

\par\bigskip
In other words:
$$ \mathcal{P}(-z) = \omega \circ \mathcal{P}(\epsilon z) \circ \omega $$
and:
$$ \mathcal{T}(-z) = \omega^* \circ \mathcal{T}(\epsilon z) \circ \omega^* $$
The following lemma tells us how the translation and multiplication operators commute:
\begin{lemma} \label{commute}
$$ \mathcal{T}(v) \mathcal{P}(w) = \Omega(vw) \mathcal{P}(w) \mathcal{T}(v) $$
\begin{proof}
\begin{align*}
\mathcal{T}(v)\mathcal{P}(w)[h(X)] & =  \mathcal{T}(v)[ \Omega(Xw) h(X) ] \\
& =  \Omega((X+v)w) h(X + v) \\
& =  \Omega(vw) \Omega(Xw) \mathcal{T}(v) [h(X)] \\
& =  \Omega(vw) \mathcal{P}(w) \mathcal{T}(v)[h(X)]
\end{align*}
\qedhere
\end{proof}
\end{lemma}

\noindent
We can now define the column operator:
\begin{equation}\label{colop}
C(z) = \sum_m (-1)^m C_m z^m = \mathcal{P}(-z) \mathcal{T}(1/z) 
\end{equation}
and the row operator:
\begin{equation}\label{rowop}
R(z) = \sum_n z^m (-1)^m R_m = \mathcal{P}(z) \mathcal{T}(-1/z) = C(-z)
\end{equation}

\begin{proposition} \label{determinant}
The Schur function has a natural expression in terms of row operators:
$$ s_\lambda(X) = R_{\lambda_1} \circ R_{\lambda_2} \circ \cdots \circ R_{\lambda_n}[1] $$
or in terms of column operators:
$$ s_\lambda(X) = C_{\lambda_1'} \circ C_{\lambda_2'} \circ \cdots \circ C_{\lambda_m'}[1] $$
\begin{proof}

Starting from the definition of the Schur function as:
$$ s_\lambda(X) = \det( e_{\lambda_i' + j - i}(X) ) $$
and expanding along the first row we get that:
\begin{align*}
s_\lambda(X) & =  \sum_{k=0}^{n-1} (-1)^{k} e_{\lambda'_1 + k}(X) s_{\mu / (k)}(X) \\
& =  \sum_{k=0}^{n-1} (-1)^k e_{\lambda'_1 + k}(X) \partial_{h_k(X)}[ s_\mu(X) ] \\
& =  C_{\lambda_1'}[s_\mu(X)] 
\end{align*}
where $\mu$ is the partition obtained from $\lambda$ by removing the first column.

\par\bigskip
Similarly, from the definition of the Schur function as:
$$ s_\lambda(X) = \det( h_{\lambda_i + j - i}(X) )$$
we get by expanding along the first row that:
$$ s_\lambda(X) = R_{\lambda_1}[s_\mu(X)] $$
where this time $\mu$ is the partition obtained from $\lambda$ by removing the first {\it row}.

\end{proof}

\end{proposition}

\begin{corollary}\label{basecase}
In the row and column cases we have:
$$s_{(n)}(X) = h_n(X)$$
$$ s_{(1^n)}(X) = e_n(X) $$
\end{corollary}

\subsubsection{Duality}

The Schur functions are, up to scaling, the unique homogeneous basis which is self dual with respect to the Hall inner product. We give here a slightly non-standard proof of this fact, which will be generalized in the next section.

\begin{proposition} \label{vertex}
The Schur functions are self-dual with respect to the Hall inner product
\begin{proof}
By proposition \ref{dual} and corollary \ref{basecase} it suffices to demonstrate that:
$$ \langle \Omega_z(Y) s_\mu(Y), s_\lambda(X) \rangle = \langle s_\mu(Y), s_\lambda(X + z) \rangle $$

\par\bigskip

\noindent
Begin by observing that:
\begin{align*}
\mathcal{T}(v) C(z) & =  \mathcal{T}(v) \mathcal{P}(-z) \mathcal{T}(1/z) \\
& =  \Omega(-vz) \mathcal{P}(-z) \mathcal{T}(v) \mathcal{T}(1/z) \\
& =  (1-vz) C(z) \mathcal{T}(v) 
\end{align*}

\noindent
which tells us that:

$$ [C(z), \mathcal{T}(v)]  = vz \, C(z) \mathcal{T}(v) $$

\par\bigskip
Equating coefficients of $(-1)^m v^k z^m$ on both sides we get that:

$$ [ C_m, \partial_{h_k} ] = - C_{m-1} \partial_{h_{k-1}} $$

\noindent
Rewriting this in the form:
$$ \partial_{h_m} C_k = C_k \partial_{h_m} + C_{k-1} \partial_{h_{m-1}} $$
we obtain:
\begin{align*}
\partial_{h_k} [s_\lambda(X)] & =  \partial_{h_k} \circ C_{\lambda_1'} \circ C_{\lambda_2'} \circ \cdots C_{\lambda_n'}[1] \\
& =  \sum_{\substack{I \subseteq [n] \\ |I| = k}} C_{\lambda_1' - a_1} \circ C_{\lambda_2' - a_2} \circ \cdots \circ C_{\lambda_n' - a_n}[1]
\end{align*}
where $a_i = 1$ if $i \in I$ and zero otherwise.

Of course it is possible that for some $k$ we have $\lambda_k = \lambda_{k+1}$ whilst $a_k \in I$ and $a_{k-1} \not \in I$, in which case the sequence of integers $(\lambda_1 - a_1, \lambda_2 - a_2, \ldots, \lambda_n - a_n)$ will not form a partition. In this case, however, we have that the determinant:
$ \det( e_{\lambda_i' - a_i + j - i}(X) ) $
contains a repeated column, and thus vanishes. 

\par\bigskip
We conclude that:
\begin{equation} \label{recurrence}
\mathcal{T}(z) [s_\lambda(X)] = s_\lambda(X + z) = \sum_{\mu \in D(\lambda)} s_\mu(X) z^{|\lambda| - |\mu|}
\end{equation}
where $D(\lambda)$ denotes all the partitions which can be obtained from $\lambda$ by removing a horizontal strip. This is the first form of the recurrence for the Schur functions.

\par\bigskip
Dually:
\begin{align*}
\mathcal{P}(v) C(z) & =  \mathcal{P}(v) \mathcal{P}(-z) \mathcal{T}(1/z) \\
& =  \mathcal{P}(-z) \mathcal{P}(v) \mathcal{T}(1/z) \\
& =  \mathcal{P}(-z) \Omega(-v/z) \mathcal{T}(1/z) \mathcal{P}(-v)  \\
& =  (1-v/z) C(z) \mathcal{P}(v) 
\end{align*}

\noindent
which tells us that:
$$ [C(z), \mathcal{P}(v)]  =  v/z \, C(z) \mathcal{P}(v) $$

\par\bigskip
By equating coefficients of $(-1)^m v^k z^m$ on both sides we find that:
$$ [ C_m, h_k(X) ] = - C_{m+1} h_{k-1}(X) $$

\noindent
Rewriting this in the form:
$$ h_m(X) C_k = C_k h_m(X) + C_{k+1} h_{m-1} $$
a similar argument to the previous reveals:
\begin{equation} \label{pieri}
\mathcal{P}(z)[s_\mu(Y)] = s_\mu(Y) \Omega_z(Y) = \sum_{\lambda \in U(\mu)} s_\lambda(Y)z^{|\lambda| - |\mu|}  
\end{equation}

where $U(\mu)$ denotes the set of all partitions which can be obtained from $\mu$ by {\it adding} a horizontal strip. 
This is the first form of the Pieri formula for the Schur functions.

\par\bigskip
Combining these two facts, we see that:
\begin{align*}
\langle s_\mu(Y) \Omega(Yz), s_\lambda(X) \rangle & =  \langle s_\mu(Y), s_\lambda(X+z)] \rangle
\end{align*}
as claimed.

\end{proof}
\end{proposition}

Note that identical arguments using the row operator rather than the column operator can be used to show that:
\begin{equation}\label{recurrence2}
\mathcal{T}(-z)[s_\lambda(X)] = s_\lambda(X - z) = \sum_{\mu \in \tilde{D}(\lambda)} s_\mu(X) (\epsilon z)^{|\lambda| - |\mu|} 
\end{equation}
where $\tilde{D}(\lambda)$ denotes all the partitions which may be obtained from $\lambda$ by removing a vertical strip.
This is the second form of the recurrence for the Schur functions.

\par\bigskip
Also we have:
\begin{equation} \label{pieri2}
\mathcal{P}(-z)[s_\mu(Y)] = s_\mu(Y) \tilde{\Omega}_z(Y) = \sum_{\lambda \in \tilde{U}(\mu)} s_\lambda(Y)(\epsilon z)^{|\lambda| - |\mu|} 
\end{equation}
where $\tilde{U}(\mu)$ denotes the set of all partitions which can be obtained from $\mu$ by adding a vertical strip. This is the second form of the Pieri formula for the Schur functions.

\par\bigskip

Putting these two facts together we see that we also have:
\begin{align*} 
\langle s_\mu(Y)\tilde{\Omega}(Yz), s_\lambda(X) \rangle & =  \langle s_\mu(Y), s_\lambda(X-z) \rangle 
\end{align*}

\par\bigskip

\subsection{Littlewood--Richardson Bases}

The Littlewood--Richardson coefficients $c^{\lambda}_{\mu \nu}$ are defined by:
$$ s_\lambda(X + X') = \sum_{\mu, \nu} c^{\lambda}_{\mu\nu} s_\mu(X) s_\nu(X')$$
These numbers may be described nicely by a combinatorial object known as {\it puzzles} \cite{KEVIN}. A very nice proof of the Littlewood-Richardson rule using ideas from Quantum Integrability can be found in \cite{ZINN}.   
In this section, we describe all the other bases $\{P_\lambda(X)\}$ of the ring of symmetric functions $\Lambda$ with the property:
$$ P_\lambda(X + X') = \sum_{\mu, \nu} c^{\lambda}_{\mu \nu}  P_\mu(X) P_\nu(X')$$
We shall say that such a basis is of {\it Littlewood--Richardson type}. 

Note that we are only interested in bases for $\Lambda$ (the ring of symmetric functions in infinitely many variables) which are stable in the sense that if $\lambda$ is a partition with at most $k$ parts, its expansion in terms of Schur functions is unchanged under the projection onto $\Lambda_k$ (the ring of symmetric functions in just $k$ variables).

\par\bigskip

\par\bigskip
We saw in the proof of proposition \ref{vertex} that:
$$ s_\lambda(X + z) = \sum_{\mu \in D(\lambda)} s_\mu(X) s_{(|\lambda| - |\mu|)}(z) $$

To prove that a given basis $\{P_\lambda(X)\}$ is of Littlewood--Richardson type it suffices to demonstrate that:
$$ P_\lambda(X + z) = \sum_{\mu \in D(\lambda)} P_\mu(X) P_{(|\lambda| - |\mu|)}(z) $$

\par\bigskip
Since the Schur functions are self-dual we have, by proposition \ref{dual}, that:
$$ s_\mu(Y) s_\nu(Y) = \sum_{\lambda} c^{\lambda}_{\mu \nu} s_\lambda(Y)$$
If $\{P_\lambda(X)\}$ is a basis of Littlewood--Richardson type, and $\{Q_\lambda(Y)\}$ is its dual basis, then we must have:
$$ Q_\mu(Y) Q_\nu(Y) = \sum_{\lambda} c^{\lambda}_{\mu \nu} Q_\lambda(Y)$$

\par\bigskip
We also saw in the proof of proposition \ref{dual} that:
$$ s_\mu(Y) 
\sum_n s_{(n)}(Y) s_{(n)}(z) 
= \sum_{\lambda \in U(\mu)} s_\lambda(Y)s_{(|\lambda| - |\mu|)}(z) $$
This will generalize to:
$$ Q_\mu(Y) 
\sum_n Q_{(n)}(Y) Q_{(n)}(z)
= \sum_{\lambda \in U(\mu)} Q_\lambda(Y)Q_{(|\lambda| - |\mu|)}(z) $$

\subsubsection{Generalized complete symmetric functions}
In the row case, the Schur functions reduce to the complete symmetric functions, which have generating function:
$$ \Omega(Xz) = \sum_n h_n(X) z^n $$
The complete symmetric functions have the property that:
$$ h_n(X + X') = \sum_k h_k(X) h_{n-k}(X') $$

\noindent
Let $f(z)$ be an arbitrary delta series, with compositional inverse $g(z)$. 

\par\bigskip
Define generalized complete symmetric functions by:
\begin{equation}
\Omega_{f}(Xz) = \prod_x \frac{1}{1-f(z)x} = \sum_n r_n(X) z^n 
\end{equation}
Also define, in the dual space:
\begin{equation}
\Phi_{g}(Yz) = \prod_y \frac{1}{1-zg(y)} = \sum_n \rho_n(Y) z^n 
\end{equation}

\par\bigskip
Since:
$$ \Omega_{f}(X + X') = \Omega_{f}(X)\Omega_{f}(X') $$
we have that:
$$ r_n(X + X') = \sum_k r_k(X) r_{n-k}(X')$$
Furthermore, we have that:
$$ r_n(X) = \sum_k \alpha_{nk} h_k(X) $$
where:
$$ \alpha_{nk} = [y^n]f(y)^k $$

\par\bigskip
In other words, $r_n(X)$ is the image of the convolution type sequence associated to $f(y)$ under co-algebra embedding of $\mathbb{Q}[x]$ into $\Lambda$ given by:
$$ \frac{x^n}{n!} \mapsto h_n(X) $$

\noindent
Note that there are plenty of other collections of ring generators $\{p_n(X)\}$ for $\Lambda$ with the property that:
$$ p_n(X + X') = \sum_k p_k(X) p_{n-k}(X') $$
For any pair of sequences $\{f_n(x)\}$ and $\{g_n(y)\}$ such:
$$ \sum_n f_n(x) g_n(y) = \frac{1}{1-xy} $$
the ring generators $\{p_n(X)\}$ defined by:
$$ \sum_n p_n(X)z^n = \prod_x f_n(x) z^n = \prod_x x^n g_n(z) $$
will have this property. However, unless they are of the special form above, it will not be possible to write them as:
$$ p_n(X) = \sum_k \Box_{nk} h_k(X) $$
instead they must have more general expansions of the form:
$$ p_n(X) = \sum_{\lambda} \Box_{\lambda} h_\lambda(X) $$
Their image under the projection: 
$$h_n(X) \mapsto \frac{x^n}{n!}$$ 
will still be a sequence of convolution type associated to some $f(y)$, however.

\par\bigskip

The idea now is to show that the same determinantal construction which produces the Schur functions from the complete symmetric functions can be used to produce more general bases with the Littlewood--Richardson property. 

In the special case where the ring bases are images of convolution type sequences, we shall see that, as a direct consequence of the multilinearity of the determinant, the associated Littlewood--Richardson bases have expansions of the form:
$$ P_\lambda(X) = \sum_{\nu \subseteq \lambda} \Box_{\lambda \nu} s_\nu(X) $$
while the dual basis will have an expansion of the form:
$$ Q_\mu(Y) = \sum_{\mu \subseteq \nu} \Box_{\nu \mu} s_\nu(Y) $$

\subsubsection{Umbral operators}
Let $U_f$ be the operator which is defined on the complete symmetric functions by:
$$ U_f[\Omega(Xz)] = \Omega_{f}(Xz)$$
extendeded multiplicatively to the whole of $\Lambda$.
Clearly we have:
$$ U_f^{-1} = U_g $$
It turns out that if $P_\lambda(X) = U_f[s_\lambda(X)]$ then $\{P_\lambda(X)\}$ is a basis of Littlewood--Richardson type with dual basis $\{Q_\lambda(Y)\}$ where $Q_\lambda(Y) = U_g^*[s_\lambda(Y)]$. 

\begin{proposition}
$$ U_g^*[\Omega(Yz)] = \Phi_{g}(Yz) $$
\begin{proof}
Let $\{n_\lambda(Y)\}$ denote the dual basis to $\{r_\lambda(X)\}$. We have:
\begin{align*}
\langle U^*_{f}[n_\mu(Y)], h_\lambda(X) \rangle 
& =  \langle n_\mu(Y), U_{f}[h_\lambda(X)] \rangle \\
& =  \langle n_\mu(Y), r_\lambda(X) \rangle \\
& =  \delta_{\lambda \mu}
\end{align*}
and so:
$$ U^*_{f}[n_\mu(Y)] = m_\mu(Y) $$
That is:
$$ n_\mu(Y) = U^*_{g}[m_\mu(Y)] $$
Now:
\begin{align*}
\Omega(XY) & =  \prod_{y} \Omega(Xy) \\
& =  \prod_{y} \sum_n r_n(X) g(y)^n \\
\end{align*}
which gives us the following explicit expression for $n_\lambda(Y)$:
$$ n_\lambda(Y) = \sum_{\sigma} \prod_k g(y_k)^{\lambda_{\sigma(k)}} $$
where the sum is over all distinct permutations of $\lambda$. Finally:
\begin{align*}
U^*_{g} [\Omega(Yz)] & =  \sum_n U^*_{g} [m_\lambda(Y)] z^{|\lambda|} \\
& =  \sum_n n_{\lambda}(Y) z^{|\lambda|} \\
& =  \prod_x \frac{1}{1-zg(y)}
\end{align*}
\qedhere
\end{proof}
\end{proposition}

\begin{proposition}
\par\bigskip
\begin{itemize}
\item $ \partial_{r_n} = U^*_g \circ \partial_{h_n} \circ U^*_f $
\item $ \partial_{\rho_n} = U_g \circ \partial_{h_n} \circ U_f $
\end{itemize}
\begin{proof}
For the first part:
\begin{align*}
\langle \partial_{r_n}[ n_\mu(Y) ], r_\lambda(X) \rangle 
& =  \langle n_\mu(Y), r_\lambda(X) r_n(X) \rangle \\
& =  \delta_{\mu, n + \lambda} \\
& =  \langle m_\mu(Y), h_\lambda(X) h_n(X) \rangle \\
& =  \langle \partial_{h_n}[m_\mu(Y)], h_\lambda(X) \rangle \\
& =  \langle \partial_{h_n} \circ U^*_{f}[n_\mu(X)], U_{g}[r_\lambda(X)] \rangle \\
& =  \langle U^*_{g} \circ \partial_{h_n} \circ U^*_{f}[n_\mu(Y)], r_\lambda(X) \rangle
\end{align*}

For the second part, let $\{\eta_\lambda(X)\}$ denote the dual basis to $\{\rho_\lambda(Y)\}$. 
As before we have:
\begin{align*}
\langle h_\mu(Y), U_{f}[\eta_\lambda(X)] \rangle 
& =  \langle U^*_{f}[h_\mu(Y)], \eta_\lambda(X),  \rangle \\
& =  \langle \rho_\mu(Y), \eta_\lambda(X)\rangle \\
& =  \delta_{\lambda \mu}
\end{align*}
and so:
$$ U_{f}[\eta_\lambda(X)] = m_\lambda(X) $$
That is:
$$ \eta_\lambda(X) = U_{g}[m_\lambda(X)] $$
Now:
\begin{align*}
\langle \rho_\mu(Y), \partial_{\rho_n}[\eta_\lambda(X)] \rangle
& =  \langle \rho_\mu(Y) \rho_{(n)}(Y), \eta_\lambda(X) \rangle \\
& =  \delta_{\mu + (n), \lambda} \\
& =  \langle h_\mu(Y) h_{(n)}(Y), m_\lambda(X) \rangle \\
& =  \langle h_\mu(Y), \partial_{h_n}[m_\lambda(X)] \rangle \\
& =  \langle U^*_g[\rho_\mu(Y)], \partial_{h_n} \circ U_f[\eta_\lambda(X)] \\
& =  \langle \rho_\mu(Y), U_g \circ \partial_{h_n} \circ U_f[\eta_\lambda(X)]
\end{align*} 
\end{proof}
\qedhere
\end{proposition}

\subsubsection{Column operators}
The (generalized) column operator for the Littlewood--Richardson basis $\{P_\lambda(X)\}$ is:
$$ C_f(z) = \sum_m (-1)^m C^f_m z^m = U_f \circ C(z) \circ U_g$$
where $C(z)$ is the column operator for the Schur functions.
\par\bigskip
The generalized column operator for the dual basis $\{Q_\lambda(Y)\}$ is:
$$ \hat{C}_f(z) = \sum_m (-1)^m \hat{C}^f_m z^m = {C}_f^*(z) = U_g^* \circ C^*(z) \circ U_f^*$$
Since these operators are adjoint, its clear that the resulting bases:
$$ P_\lambda(X) = {C}^f_{\lambda_1'} \circ {C}^f_{\lambda_2'} \circ \cdots \circ {C}^f_{\lambda_n'}[1]$$
$$ Q_\lambda(Y) = \hat{C}^f_{\lambda_1'} \circ \hat{C}^f_{\lambda_2'} \circ  \cdots \circ \hat{C}^f_{\lambda_n'}[1]$$
are dual. More explicitly:
\begin{align*}
\langle Q_\mu(Y), P_\lambda(X) \rangle & =  \langle U^*_g[s_\mu(Y)], U_f[s_\lambda(X)] \rangle \\
& =  \langle s_\mu(Y), U_g \circ U_f[s_\lambda(X) \rangle \\
& =  \langle s_\mu(Y), s_\lambda(X) \rangle \\
& =  \delta_{\lambda \mu}
\end{align*}

\par\bigskip
Note that by the same argument used in proposition \ref{determinant} we have:
$$ P_\lambda(X) = \det(r_{\lambda_i + j - i}(X))$$
$$ Q_\lambda(Y) = \det(\rho_{\lambda_i + j - i}(Y))$$

\noindent
Now define the generalized multiplication operator:
\begin{align*}
\mathcal{P}_{f}(z) & =  \sum_k z^k r_k(X) \\
& =  \sum_k f(z)^k h_k(X) \\
& =  U_f \circ \mathcal{P}(z) \circ U_g 
\end{align*}
as well as the generalized translation operator:
\begin{align*}
\mathcal{T}_{g}(z) & =  \sum_k z^k \partial_{\rho_k} \\
& =  \sum_k r_k(z) \partial_{h_k} \\
& =  U_f \circ \mathcal{T}(z) \circ U_g
\end{align*}

\par\bigskip
Also let us define, in the dual space:
\begin{align*}
\hat{\mathcal{P}}_{g}(z) & =  \mathcal{T}_g^*(z) \\
& =  U_g^* \circ \mathcal{P}(z) \circ U_f^* \\
& =  \sum_k z^k \rho_k(Y) \\
& =  \sum_k r_k(z) h_k(X)
\end{align*}
as well as:
\begin{align*}
\hat{\mathcal{T}}_{f}(z) & =  \mathcal{P}_f^*(z) \\
& =  U^*_g \circ \mathcal{T}(z) \circ U^*_f \\
& =  \sum_k z^k \partial_{r_k} \\
& =  \sum_k f(z)^k \partial_{h_k}
\end{align*}

\par\bigskip
We have now that:
$$ {C}_f(z) = \mathcal{P}_{-f}(z) \circ \mathcal{T}_{1/g}(z)$$
while:
$$ \hat{C}_f(z) = \hat{\mathcal{P}}_{-f}(z) \circ \hat{\mathcal{T}}_{1/g}(z) $$

\begin{proposition}
The (non-homogeneous) basis for the ring of symmetric functions defined by:
$$ P_\lambda(X) = C^f_{\lambda_1'} \circ C^f_{\lambda_2'} \circ \cdots \circ C^f_{\lambda_m'}[1]$$
is of Littlewood--Richardson type.
\begin{proof}
We must show that:
$$ \langle Q_\lambda(Y) \Phi_g(Yz), P_\mu(X) \rangle = \langle Q_\lambda(Y), P_\mu(X+z) \rangle $$
By conjugating Lemma \ref{commute} by $U_f$ on the left and $U_g$ on the right, we find that:
$$ \mathcal{T}_{g}(v) \circ \mathcal{P}_{f}(w) = \Omega(vw) \circ \mathcal{P}_{f}(w) \circ \mathcal{T}_{g}(v) $$
This in turn implies that:
$$ [ \tilde{C}(z), \mathcal{T}_{g}(v)] = vz \circ \tilde{C}(z) \circ \mathcal{T}_{g}(v) $$
By arguments which are essentially identical to those in the proof of proposition \ref{vertex} we obtain the following generalized recurrence:
$$ \partial_{\rho_k}[P_\lambda(X)] = \sum_{\mu \in D_k(\lambda)} P_\mu(X)$$
Next by using the fact that:
$$ \mathcal{T}(z) = \sum_k r_k(z) \partial_{\rho_k} $$
we can re-write the above as:
$$ \mathcal{T}(z)[P_\lambda(X)] = \sum_{\mu \in D(\lambda)} P_\mu(X) P_{(|\lambda| - |\mu|)}(z) $$
\par\bigskip
Similarly, by conjugating Lemma \ref{commute} on the left by $U^*_g$ and the right by $U^*_f$ we find that:
$$ \hat{\mathcal{T}}_{g}(v) \circ \hat{\mathcal{P}}_{f}(w) = \Omega(vw) \circ \hat{\mathcal{P}}_{f}(w) \circ \hat{\mathcal{T}}_{g}(v) $$
This tells us that:
$$ [  \hat{C}(z), \hat{\mathcal{P}}_{f}(v) ] = v/z \circ \hat{C}(z) \circ \hat{\mathcal{P}}_{f}(v) $$
and so we have the following Pieri formula	:
$$ Q_\mu(Y) \rho_k(Y)  = \sum_{\lambda \in U_k(\mu)} Q_\lambda(Y)$$
That is:
$$ \hat{\mathcal{P}}_{f}(z)[Q_\mu(Y)] = Q_\mu(Y) \Phi_g(Yz) = \sum_{\lambda \in U(\mu)} Q_\lambda(Y) z^{|\lambda| - |\mu|	} $$
Putting these two facts together we find that:
$$ \langle Q_\lambda(Y) \Phi_g(Yz), P_\mu(X) \rangle = \langle Q_\lambda(Y), P_\mu(X+z) \rangle $$
as claimed.
\end{proof}
\end{proposition}

\subsubsection{Generalized elementary symmetric functions}
Now define generalized elementary symmetric functions:
$$ \tilde{\Omega}_{f}(Xz) = \prod_x (1 - f(-z)x) = \sum_n c_n(X) z^n $$
and in the dual space:
$$ \tilde{\Phi}_{g}(Yz) = \prod_y (1+zg(y)) = \sum_n \gamma_n(Y) z^n $$
A minor modification of the argument at the end of section \ref{schurfun} can be used to show that:
$$ \det(r_{\lambda_i + j - i}(X)) = \det(c_{\lambda'_i + j - i}(\epsilon X)) $$
$$ \det(\rho_{\lambda_i + j - i}(Y)) = \det(\gamma_{\lambda'_i + j - i}(\epsilon Y)) $$
\par\bigskip
The generalized row operators are: 
$${R}_f(z) = {C}_f(-z)$$ 
$$\hat{R}_f(z) = \hat{C}_f(-z)$$ 
One can show that:
$$ P_\lambda(X) = {R}^f_{\lambda_1} \circ {R}^f_{\lambda_2} \circ \cdots \circ {R}^f_{\lambda_n}[1] $$
$$ Q_\lambda(Y) = \hat{R}^f_{\lambda_1} \circ \hat{R}^f_{\lambda_2} \circ \cdots \circ \hat{R}^f_{\lambda_n}[1] $$
which implies that:
$$ \mathcal{T}(-z)[P_\lambda(X)] = P_\lambda(X - z) = \sum_{\mu \in \tilde{D}_{\lambda}} P_\mu(X) P_{(|\lambda| - |\mu|)}(\epsilon z) $$
$$ \hat{\mathcal{P}}_{-f}(z)[Q_\mu(Y)] = Q_\mu(Y) \tilde{\Phi}_f(Yz) = \sum_{\lambda \in \tilde{U}(\mu)} Q_\lambda(Y) Q_{(|\lambda| - |\mu|)}(\epsilon z) $$
In other words:
$$ \langle Q_\lambda(Y)\tilde{\Phi}_f(Yz) , P_\mu(X) \rangle = \langle Q_\lambda(Y), P_\mu(X-z) \rangle $$

\subsection{Examples}

In this section we shall give symmetric function analogues of the rising and falling factorials, as well as the rook polynomials, leading to a generalized family of Stirling and Lah numbers which are indexed by partitions rather than integers. 

\par\bigskip
For any formal power series $f(z)$ let $\bar{f}(z) = f(-z)$. 


\par\bigskip
Let us use the following notation for the falling factorial:
$$ (x)_n = x(x-1)(x-2) \cdots (x-n+1)$$
and the rising factorial:
$$ (x)^n = x(x+1)(x+2) \cdots (x+n-1) $$
The Stirling numbers of the first kind are defined by:
$$ (x)_n = \sum_k s(n,k) x^k $$
We have:
$$ (x)^n = \sum_k |s(n,k)| x^k $$
One may check that $\{\frac{(x)_n}{n!}\}$ is the sequence of convolution type associated to $f(z) = (\exp(z) - 1)$ while $\{\frac{(x)^n}{n!}\}$ is the sequence of convolution type associated to $-\bar{f}(z) = (1 - \exp(-z))$.




\bigskip
Now let $\{A_\lambda(X)\}$ denote the Littlewood--Richardson basis for the ring of symmetric functions $\Lambda$ that is associated to $f(z) = \exp(z) - 1$. 
The generating function for the row case is:
$$ \Omega_f(Xz) = \prod_{x \in X} \frac{1}{1 - x \log (1+z)}  $$
while the generating function for the column case is:
$$ \tilde{\Omega}_f(Xz) = \prod_{x \in X} (1 - x \log(1-z)) $$

Similarly let $\{B_\lambda(X)\}$ denote the Littlewood--Richardson basis associated to $-\bar{f}(z) = 1 - \exp(-z)$. The generating function for the row case is:
\begin{align*}
\Omega_{-\bar{f}}(Xz)
& =  \prod_{x \in X} \frac{1}{1 - x  \log \left ( \frac{1}{1-z} \right )}  \\
& =  \prod_{x \in X} \frac{1}{1 + x \log(1-z)} \\
& =  \tilde{\Omega}_f(-\epsilon Xz) \\
& =  \omega[ \tilde{\Omega}_f(Xz)]
\end{align*}
The generating function for the column case is:
\begin{align*}
\tilde{\Omega}_{-\bar{f}}(Xz) 
& =  \prod_{x \in X} (1 - x \log \left ( \frac{1}{1+z} \right )) \\ 
& =  \prod_{x \in X} (1 + x \log (1 + z)) \\
& =  \Omega_f(- \epsilon Xz) \\
& =  \omega[ \Omega_f(Xz) ]
\end{align*}

Let $A = (a_{\lambda \mu})$ denote the transition matrix from the basis $\{A_{\lambda}(X)\}$ to the usual Schur basis $\{s_\lambda(X)\}$, and let $B = (b_{\lambda \mu})$ denote the transition matrix from $\{B_\lambda(X)\}$ to $\{s_\lambda(X)\}$. On the next page we give the corner of these matrices corresponding to partitions with at most five parts.

\par\bigskip
The order which we are using on partitions is:
\begin{gather*}
\{ \{1\} \}
\\
\{ \{2\}, \{1, 1\} \}
\\
\{ \{3\}, \{2, 1\}, \{1, 1, 1\} \}
\\
 \{ \{4\}, \{3, 1\}, \{2, 2\}, \{2, 1, 1\}, \{1, 1, 1, 1\} \} 
\\
\{ \{5\}, \{4, 1\}, \{3, 2\}, \{3, 1, 1\}, \{2, 2, 1\}, \{2, 1, 1, 1\}, \{1, 1, 1, 1, 1\} \}
\end{gather*}

\bigskip
\noindent
Here is a portion of the matrix $A$:
\[
\left(
\setlength{\arraycolsep}{2pt}
\def\str{{\vrule height12pt depth5pt width0pt}}
\begin{array}{c|cc|ccc|ccccc|ccccccc}
 1 & -\frac{1}{2} & \frac{1}{2} & \frac{1}{3} & -\frac{1}{3} & \frac{1}{3} & -\frac{1}{4} & \frac{1}{4} & 0 & -\frac{1}{4} & \frac{1}{4} &   \frac{1}{5} & -\frac{1}{5} & 0 & \frac{1}{5} & 0 & -\frac{1}{5} & \frac{1}{5} \\
[2pt]\hline\str 
 0 & 1 & 0 & -1 & \frac{1}{2} & 0 & \frac{11}{12} & -\frac{7}{12} & -\frac{1}{12} & \frac{1}{3} & 0 & -\frac{5}{6} & \frac{7}{12} &   \frac{1}{12} & -\frac{5}{12} & -\frac{1}{12} & \frac{1}{4} & 0 \\
 \vdots & \ddots & 1 & 0 & -\frac{1}{2} & 1 & 0 & \frac{1}{3} & -\frac{1}{12} & -\frac{7}{12} & \frac{11}{12} & 0 & -\frac{1}{4} & \frac{1}{12} &  \frac{5}{12} & -\frac{1}{12} & -\frac{7}{12} & \frac{5}{6} \\
[2pt]\hline\str 
   &  &   & 1 & 0 & 0 & -\frac{3}{2} & \frac{1}{2} & 0 & 0 & 0 & \frac{7}{4} & -\frac{5}{6} & -\frac{1}{12} & \frac{1}{3} & 0 & 0 & 0 \\
   &   &   &   & 1 & 0 & 0 & -1 & 0 & 1 & 0 & 0 & \frac{11}{12} & -\frac{1}{12} & -\frac{13}{12} & -\frac{1}{12} & \frac{11}{12} & 0 \\
   &   &   &   &   & 1 & 0 & 0 & 0 & -\frac{1}{2} & \frac{3}{2} & 0 & 0 & 0 & \frac{1}{3} & -\frac{1}{12} & -\frac{5}{6} & \frac{7}{4} \\
[2pt]\hline\str 
   &   &   &   &   &   & 1 & 0 & 0 & 0 & 0 & -2 & \frac{1}{2} & 0 & 0 & 0 & 0 & 0 \\
   &   &   &   &   &   &   & 1 & 0 & 0 & 0 & 0 & -\frac{3}{2} & 0 & 1 & 0 & 0 & 0 \\
   &   &   &   &   &   &   &   & 1 & 0 & 0 & 0 & 0 & -1 & 0 & 1 & 0 & 0 \\
   &   &   &   &   &   &   &   &   & 1 & 0 & 0 & 0 & 0 & -1 & 0 & \frac{3}{2} & 0 \\
   &   &   &   &   &   &   &   &   &   & 1 & 0 & 0 & 0 & 0 & 0 & -\frac{1}{2} & 2 \\
[2pt]\hline\str 
   &   &   &   &   &   &   &   &   &   &   & 1 & 0 & 0 & 0 & 0 & 0 & 0 \\
   &   &   &   &   &   &   &   &   &   &   &   & 1 & 0 & 0 & 0 & 0 & 0 \\
   &   &   &   &   &   &   &   &   &   &   &   &   & 1 & 0 & 0 & 0 & 0 \\
   &   &   &   &   &   &   &   &   &   &   &   &   &   & 1 & 0 & 0 & 0 \\
   &   &   &   &   &   &   &   &   &   &   &   &   &   &   & 1 & 0 & 0 \\
 \vdots &   &   &   &   &   &   &   &   &   &   &   &   &   &   & \ddots  & 1 & 0 \\
 0 & \cdots &   &   &   &   &   &   &   &   &   &   &   &   &   & \cdots & 0 & 1
\end{array}
\right)
\]
The columns are the images of the Schur functions under the transformation $h_k \mapsto \sum_{i \leq k} c_{k,i} h_i$ where $c_{k,i}$ are the coefficients of the powers of $\log(1+z)$.

\par\bigskip
Here is a portion of the matrix $B$:
\[
\left(
\setlength{\arraycolsep}{2pt}
\def\str{{\vrule height12pt depth5pt width0pt}}
\begin{array}{c|cc|ccc|ccccc|ccccccc}
 1 & \frac{1}{2} & -\frac{1}{2} & \frac{1}{3} & -\frac{1}{3} & \frac{1}{3} & \frac{1}{4} & -\frac{1}{4} & 0 & \frac{1}{4} & -\frac{1}{4} &   \frac{1}{5} & -\frac{1}{5} & 0 & \frac{1}{5} & 0 & -\frac{1}{5} & \frac{1}{5} \\
[2pt]\hline\str
 0 & 1 & 0 & 1 & -\frac{1}{2} & 0 & \frac{11}{12} & -\frac{7}{12} & -\frac{1}{12} & \frac{1}{3} & 0 & \frac{5}{6} & -\frac{7}{12} &   -\frac{1}{12} & \frac{5}{12} & \frac{1}{12} & -\frac{1}{4} & 0 \\
 \vdots & \ddots & 1 & 0 & \frac{1}{2} & -1 & 0 & \frac{1}{3} & -\frac{1}{12} & -\frac{7}{12} & \frac{11}{12} & 0 & \frac{1}{4} & -\frac{1}{12} &   -\frac{5}{12} & \frac{1}{12} & \frac{7}{12} & -\frac{5}{6} \\
[2pt]\hline\str
   &   &   & 1 & 0 & 0 & \frac{3}{2} & -\frac{1}{2} & 0 & 0 & 0 & \frac{7}{4} & -\frac{5}{6} & -\frac{1}{12} & \frac{1}{3} & 0 & 0 & 0 \\
   &   &   &   & 1 & 0 & 0 & 1 & 0 & -1 & 0 & 0 & \frac{11}{12} & -\frac{1}{12} & -\frac{13}{12} & -\frac{1}{12} & \frac{11}{12} & 0 \\
   &   &   &   &   & 1 & 0 & 0 & 0 & \frac{1}{2} & -\frac{3}{2} & 0 & 0 & 0 & \frac{1}{3} & -\frac{1}{12} & -\frac{5}{6} & \frac{7}{4} \\
[2pt]\hline\str
   &   &   &   &   &   & 1 & 0 & 0 & 0 & 0 & 2 & -\frac{1}{2} & 0 & 0 & 0 & 0 & 0 \\
   &   &   &   &   &   &   & 1 & 0 & 0 & 0 & 0 & \frac{3}{2} & 0 & -1 & 0 & 0 & 0 \\
   &   &   &   &   &   &   &   & 1 & 0 & 0 & 0 & 0 & 1 & 0 & -1 & 0 & 0 \\
   &   &   &   &   &   &   &   &   & 1 & 0 & 0 & 0 & 0 & 1 & 0 & -\frac{3}{2} & 0 \\
   &   &   &   &   &   &   &   &   &   & 1 & 0 & 0 & 0 & 0 & 0 & \frac{1}{2} & -2 \\
[2pt]\hline\str
   &   &   &   &   &   &   &   &   &   &   & 1 & 0 & 0 & 0 & 0 & 0 & 0 \\
   &   &   &   &   &   &   &   &   &   &   &   & 1 & 0 & 0 & 0 & 0 & 0 \\
   &   &   &   &   &   &   &   &   &   &   &   &   & 1 & 0 & 0 & 0 & 0 \\
   &   &   &   &   &   &   &   &   &   &   &   &   &   & 1 & 0 & 0 & 0 \\
   &   &   &   &   &   &   &   &   &   &   &   &   &   &   & 1 & 0 & 0 \\
 \vdots  &   &   &   &   &   &   &   &   &   &   &   &   &   &   & \ddots  & 1 & 0 \\
 0  & \cdots &   &   &   &   &   &   &   &   &   &   &   &   &   & \cdots & 0 & 1
\end{array}
\right)
\]
A similar remark applies to $B$ with $\log(1+z)$ replaced with $-\log(1-z)$.
Observe that $b_{\lambda \mu} = a_{\lambda'\, \mu'}$, 
and 
$b_{\lambda \mu} = (-1)^{|\lambda| - |\mu|} a_{\lambda  \mu}$.

If one looks at just entries of $A$ corresponding to row diagrams (scaled by $n! / k!$), or alternatively those of $B$ corresponding to column diagrams, one recovers the Stirling numbers of the first kind:
$$ \left(
\begin{array}{rrrrr}
 1 & -1 & 2 & -6 & 24 \\
 0 & 1 & -3 & 11 & -50 \\
 0 & 0 & 1 & -6 & 35 \\
 0 & 0 & 0 & 1 & -10 \\
 0 & 0 & 0 & 0 & 1
\end{array}
\right) $$
Switching the roles of $A$ and $B$,
one obtains the unsigned Stirling numbers of the first kind:
$$ \left(
\begin{array}{rrrrr}
 1 & 1 & 2 & 6 & 24 \\
 0 & 1 & 3 & 11 & 50 \\
 0 & 0 & 1 & 6 & 35 \\
 0 & 0 & 0 & 1 & 10 \\
 0 & 0 & 0 & 0 & 1
\end{array}
\right) $$

\par\bigskip

\bigskip
\noindent
Let us define the Lah numbers by:
$$ (x)_n = \sum_k L(n,k) (x)^k $$

\par\bigskip
Let $h(z) = \frac{z}{1+z}$ and observe that $h^{-1}(z) = \frac{z}{1-z}$ 

We have:
\begin{align*}
-\bar{g}(h(z))   &=  \log \left (\frac{1}{1 - \frac{z}{z+1}} \right) \\
& =  \log(1+z) \\
& =  g(z) \\
\end{align*}
In other words:
$$ U_{-\bar{f}} \circ U_{h^{-1}} = U_{f} $$
Suppose now that the rook polynomials have the expansion:
$$ l_n(x) = U_{h^{-1}} \left [ \frac{x^n}{n!} \right ] = \sum_k a_{nk} \frac{x^k}{k!} $$
then we have:
$$ U_{-\bar{f}} \left[ l_n(x) \right ] = \sum_k a_{nk} \frac{(x)^k}{k!} $$
while:
$$ U_{f} \left [ \frac{x^n}{n!} \right ] = \frac{(x)_n}{n!} $$
Thus we must have:
$$l_n(x) = \sum_k |L(n,k)| x^k $$

\bigskip
\noindent
Now consider the Littlewood--Richardson basis $\{L_\lambda(X)\}$ associated to $h(z)=\frac{z}{1-z}$. Here is a portion of the transition matrix $L = (l_{\lambda \mu})$:
\[
\left(
\setlength{\arraycolsep}{2pt}
\def\str{{\vrule height12pt depth5pt width0pt}}
\begin{array}{c|cc|ccc|ccccc|ccccccc}
 1 & 1 & -1 & 1 & -1 & 1 & 1 & -1 & 0 & 1 & -1   & 1 & -1 & 0 & 1 & 0 & -1 & 1 \\
[2pt]\hline\str
 0 & 1 & 0 & 2 & -1 & 0 & 3 & -2 & 0 & 1 & 0 &   4 & -3 & 0 & 2 & 0 & -1 & 0 \\
 \vdots & \ddots & 1 & 0 & 1 & -2 & 0 & 1 & 0 & -2 & 3 & 0   & 1 & 0 & -2 & 0 & 3 & -4 \\
[2pt]\hline\str
  &  &  & 1 & 0 & 0 & 3 & -1 & 0 & 0 & 0 & 6 &   -3 & 0 & 1 & 0 & 0 & 0 \\
  &  &  &  & 1 & 0 & 0 & 2 & 0 & -2 & 0 & 0 &   3 & 0 & -4 & 0 & 3 & 0 \\
  &  &  &  &  & 1 & 0 & 0 & 0 & 1 & -3 & 0 & 0   & 0 & 1 & 0 & -3 & 6 \\
[2pt]\hline\str
  &  &  &  &  &  & 1 & 0 & 0 & 0 & 0 & 4 & -1   & 0 & 0 & 0 & 0 & 0 \\
  &  &  &  &  &  &  & 1 & 0 & 0 & 0 & 0 & 3 &   0 & -2 & 0 & 0 & 0 \\
  &  &  &  &  &  &  &  & 1 & 0 & 0 & 0 & 0 & 2   & 0 & -2 & 0 & 0 \\
  &  &  &  &  &  &  &  &  & 1 & 0 & 0 & 0 & 0   & 2 & 0 & -3 & 0 \\
  &  &  &  &  &  &  &  &  &  & 1 & 0 & 0 & 0 &   0 & 0 & 1 & -4 \\
[2pt]\hline\str
  &  &  &  &  &  &  &  &  &  &  & 1 & 0 & 0 &   0 & 0 & 0 & 0 \\
  &  &  &  &  &  &  &  &  &  &  &  & 1 & 0 & 0   & 0 & 0 & 0 \\
  &  &  &  &  &  &  &  &  &  &  &  &  & 1 & 0   & 0 & 0 & 0 \\
  &  &  &  &  &  &  &  &  &  &  &  &  &  & 1 &   0 & 0 & 0 \\
  &  &  &  &  &  &  &  &  &  &  &  &  &  &  &   1 & 0 & 0 \\
 \vdots &  &  &  &  &  &  &  &  &  &  &  &  &  &  & \ddots  & 1 & 0 \\
 0 & \cdots &  &  &  &  &  &  &  &  &  &  &  &  &  & \cdots  & 0 & 1
\end{array}
\right)
\]
The columns are the images under the Schur functions under the transformation $h_k \mapsto \sum_{i \leq k} {k\choose i}   h_i$.

\par\bigskip
Looking at just the rows (scaled by $n! / k!$), we recover the Lah numbers:
$$ \left(
\begin{array}{rrrrr}
 1 & 2 & 6 & 24 & 120 \\
 0 & 1 & 6 & 36 & 240 \\
 0 & 0 & 1 & 12 & 120 \\
 0 & 0 & 0 & 1 & 20 \\
 0 & 0 & 0 & 0 & 1
\end{array}
\right) $$

\par\bigskip
One may check that $AL = B$.
\newpage
\section{A generating function identity for Macdonald polynomials}

The Macdonald polynomials are a $qt$-deformation of the Schur functions. The main result of the second part of this thesis is the proof of a generating function identity for Macdonald polynomials which was originally conjectured by Kawanaka \cite{KAWANAKA}. 

In section \ref{intro} we define the Macdonald polynomials and state the Kawanaka conjecture. Section \ref{resultants} contains some technical lemmas which are not needed until the final step of the proof. In section \ref{pierirecurrence} we discuss the Pieri formula and recurrence for the Macdonald formulas, which generalize those for the Schur functions which were discussed in detail in the first part of this thesis. The proof itself is contained in section \ref{proof}.

\subsection{Macdonald Polynomials}\label{intro}
\subsubsection{Notation}
Recall from section \ref{pleth} the generating series for the complete symmetric functions:
\[
\Omega(zX)=\prod_{x\in X}\frac{1}{1-z\,x}
\]

Here we again use the plethystic notation. 
Inside a symmetric function the expression $\frac{1-t}{1-q}$ denotes the alphabet:
$$ \frac{1 - t}{1 - q} = \{1, q, q^2, \ldots -t, -tq, -tq^2, \ldots \} $$	

Recall also the important distinction between plethystic negation:
\begin{equation*}
\Omega(-X) = \prod_{x \in X} (1 - x) 
\end{equation*}
and formal negation:
\begin{equation*}
\Omega(\epsilon X) = \prod_{x \in X} \frac{1}{1 + x}
\end{equation*}

We remark that the operator $\Omega$ acting on alphabets plays a role in symmetric function theory that is in many ways analogous to that of the exponential function in the theory of functions of a single variable. In particular:
$$ \Omega(X + Y) = \Omega(X) \Omega(Y) $$
Also note that:
$$ \Omega \left ( \frac{z}{1-q} \right ) = \frac{1}{(z;q)_\infty} = \sum_k \frac{z^k}{(q;q)_k} $$
is the first $q$-exponential function, mentioned in section \ref{partitions}, whilst:
$$ \Omega \left ( \frac{-z}{1-q} \right ) = (-z;q)_\infty = \sum_k q^{\binom{k}{2}} \frac{z^k}{(q;q)_k}$$
is the second $q$-exponential function.

\subsubsection{Operator definition}

\noindent
In this section we shall be workng with a finite alphabet $X_n$ with exactly $n$ letters:
$$ X_n = x_1 + x_2 + \cdots + x_n $$
Let $\Delta(X_n)$ denote the Vandemonde determinant:

$$ \Delta(X_n) = \left \|
\begin{array}{cccc}
1 & 1 & \ldots & 1 \\
x_1 & x_2 & \ldots & x_n \\
\ldots & \ldots & \ldots & \ldots \\
x_1^{n-1} & x_2^{n-1} & \ldots & x_n^{n-1}
\end{array}
\right \| = \prod_{i < j} (x_j - x_i) $$

\medskip
\noindent
The Macdonald polynomials $\{ P_\lambda(X_n) \}$ were originally introduced in \cite{MAC} as the eigenfunctions of the operator $D$:
$$ f(x_1, x_2, \ldots x_n) \mapsto \sum_{i = 1}^n \prod_{j \neq i} \frac{(tx_i - x_j)}{(x_i - x_j)} f(x_1, \ldots, q x_i, \ldots x_n) $$
with eigenvalues:
$$ \sum_{i=1}^n q^{\lambda_i} t^{n-i} $$
We may rewrite this operator using the plethystic notation \cite{LASCOUX, GARSIA1, GARSIA2} as:
$$ f(X_n) \mapsto \frac{1}{\Delta(X_n)} \sum_{i=1}^n \Delta(X_n + (t-1)x_i) f(X_n + (q-1)x_i) $$
Although it is not possible to take the limit $n \to \infty$ without the eigenvalues diverging, by using a slightly modified version of this operator, namely:
$$ f(X_n) \mapsto t^{-n}  \frac{1}{\Delta(X_n)} \sum_{i=1}^n  \Bigl ( \Delta(X_n + (t-1)x_i) f(X_n + (q-1)x_i) \Bigr ) - \sum_{i=1}^n t^{-i} f(X_n) $$
which has eigenvalues
$ \sum_{i = 1}^n ( q^{\lambda_i} - 1)t^{-i} $,
Macdonald was able to show \cite{BLACK} (page 321) that his polynomials have the property that:
$$ P_\lambda(X_n) = \pi_n[ P_\lambda(X_{n+1}) ] $$
where $\pi_n : \Lambda_{n+1} \to \Lambda_n$ is the map which sets the $(n+1)$th variable equal to zero.

\subsubsection{Characterization using the inner product}
\noindent
Macdonald showed in \cite{MAC} that his operator $D$ is self-adjoint (and thus its eigenfunctions are orthogonal) with respect to the deformed Hall inner product:
$$ \langle f(Y), g(X) \rangle_{q,t} = \left \langle f \left( Y \frac{1-t}{1-q} \right ), g(X) \right \rangle $$
Dual bases $\{P_\lambda(X)\}$ and $\{Q_\lambda(Y)\}$ with respect to the Macdonald inner product are characterized by the property that:
$$ \sum_\lambda P_\lambda(X) Q_\lambda(Y) = \Omega \left(XY\frac{1-t}{1-q} \right) = \prod_{\substack{x \in X \\ y \in Y}} \frac{(txy,q)_\infty}{(xy,q)_\infty} $$

Note the important fact that:
$$ \Omega \left ( (X + Y) \frac{1-t}{1-q} \right ) = \Omega \left (X \frac{1-t}{1-q} \right ) \Omega \left ( Y \frac{1-t}{1-q} \right ) $$

The Macdonald polynomials may be characterized as the unique basis for the ring of symmetric functions over $\mathbb{Q}(q,t)$ which is both orthogonal (but not orthonormal) with respect to the Macdonald inner product, and whose expansion in terms of the monomial symmetric functions is strictly upper triangular with respect to the dominance order \cite{MAC} on partitions:
$$ P_\lambda(X) = \sum_{\mu \triangleleft \lambda} \Box_{\lambda \mu} m_\mu(X) $$




In particular in the column case the Macdonald polynomials correspond to the elementary symmetric functions:
$$P_{(1^n)}(X) = e_n(X)$$

\subsubsection{Arms and legs}
For $s \in \lambda$ some box in a partition $\lambda$ the {\it arm length} $a_\lambda(s)$ is defined to be the number of boxes in $\lambda$ lying directly to the right of the box $s$, while the {\it leg length} is defined to be the number of boxes in the partition $\lambda$ lying directly below the box $s$.

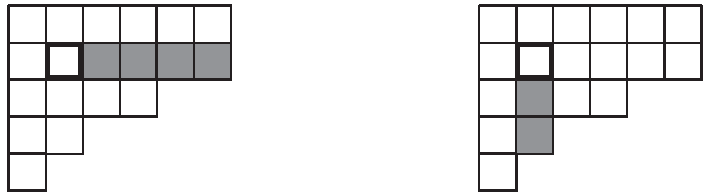
\[
\tableau{ \ &\ &\ &\ &\ &\ \\ \ &\thickcell &\graycell &\graycell &\graycell &\graycell \\ \ &\ &\ &\ \\ \ &\ \\ \ }
\qquad \qquad \qquad \qquad
\tableau{ \ &\ &\ &\ &\ &\ \\ \ &\thickcell & \ &\ &\ &\ \\ \ &\graycell &\ &\ \\ \ &\graycell \\ \ }
\]

\noindent
In other words, if $s = (i,j)$ then $a_\lambda(s) = \lambda_i - j$ and $l_\lambda(s) = \lambda'_j - i$. If the box $s$ lies outside the partition $\lambda$ then we defined $a_\lambda(s) = l_\lambda(s) = 0$.

\medskip
\noindent
Let $n(\lambda) = \sum_{s \in \lambda} a_\lambda(s)$ and let $\tilde{n}(\lambda) = \sum_{s \in \lambda} l_\lambda(s) = n(\lambda')$.

\medskip
\noindent
Let us define:
$$ B_\lambda(q,t) = \sum_\lambda q^{a_\lambda(s)} t^{l_\lambda(s)} $$
We have, of course, that:
$$B_{\lambda'}(t,q) = B_\lambda(q,t)$$

\medskip
\noindent
It is a surprising fact \cite{BLACK} (pages 338,339) that:
$$ \langle P_\lambda(Y), P_\lambda(X) \rangle_{q,t} = \frac{\prod_{s \in \lambda}(1 - q^{a_\lambda(s)+1} t^{l_\lambda(s)})}{\prod_{s \in \lambda}(1 - q^{a_\lambda(s)} t^{l_\lambda(s)+1})} 
= \Omega ( (t-q) B_\lambda(q,t) ) $$

\subsubsection{Duality}
The Macdonald $Q$-functions are defined to be dual to the Macdonald $P$-functions. Since the Macdonald $P$-functions are orthogonal with respect to the Hall inner product we have:
$$ Q_\lambda(X) = \frac{P_\lambda(X)}{\langle P_\lambda(Y) P_\lambda(X) \rangle_{q,t}} $$
The operator $\omega_{q,t}$ given by:
$$ \omega_{q,t} [f(X)] = \omega f\left(\frac{1-q}{1-t}X\right) $$
is self-adjoint with respect to the Macdonald inner product.

\medskip
\noindent
It is another surprising fact  \cite{BLACK} (page 327) that:
$$ \omega_{q,t}[ P_\lambda(X;q,t) ] = Q_{\lambda'}(X;t,q) $$

Using this, one can show that in the one row case, the Macdonald polynomials correspond to the modified version of the complete symmetric functions \cite{BLACK} (page 311) defined by:
$$P_{(n)}(X) = \frac{(q;q)_n}{(t;q)_n} g_n(X)$$
Where:
$$ \sum_n g_n(X)z^n = \Omega_z \left( X\frac{1-t}{1-q} \right)$$

\subsubsection{Kawanaka conjecture}
In part II of this thesis we shall prove the following generating function identity for the Macdonald polynomials, which was originally conjectured by Kawanaka \cite{KAWANAKA} and proved in the case of Hall-Littlewood Polynomials $(q=0)$:
$$ \sum_\lambda  \left ( 
\prod_{s \in \lambda} \frac{1 + q^{a_\lambda(s)} t^{l_\lambda(s)+1}}{1 - q^{a_\lambda(s)+1} t^{l_\lambda(s)}} \right )
P_\lambda(X;q^2,t^2) = \prod_i \frac{(-t x_i ; q)_{\infty}}{(x_i;q)_{\infty}} \prod_{i < j} \frac{(t^2 x_i x_j; q^2)_\infty}{(x_i x_j; q^2)_\infty}$$

\medskip
\noindent
By $P_\lambda(X;q^2,t^2)$ we mean the Macdonald polynomial $P_\lambda(X)$ for which every occurance of the variable $q$ has been replaced by $q^2$ and similarly every occurence of the variable $t$ has been replaced by $t^2$. 

\medskip
This identity complements the following two generating function identities for Macdonald polynomials which can be found in the case where $b=0$ or $b=1$ b on page 349 of Macdonald \cite{BLACK}:

$$ \sum_\lambda b^{c(\lambda)} \prod_{\substack{s \in \lambda \\ l_\lambda(s) \,\, {\text even}}}\frac{1- q^{a_\lambda(s)}t^{l_\lambda(s) + 1}}{1 - q^{a_\lambda(s) + 1} t^{l_\lambda(s)}}P_\lambda(X;q,t) = \prod_{i \leq j} \prod_i \frac{(btx_i;q)_\infty}{(bx_i)_\infty} \frac{(t x_i x_j ; q)_\infty}{(x_i x_j; q)_\infty} $$
$$ \sum_\lambda b^{r(\lambda)} \prod_{\substack{s \in \lambda \\ l_\lambda(s) \,\, {\text even}}}\frac{1- q^{a_\lambda(s)}t^{l_\lambda(s) + 1}}{1 - q^{a_\lambda(s) + 1} t^{l_\lambda(s)}} P_\lambda(X;q,t) = \prod_{i \leq j} \prod_i \frac{(btx_i;q)_\infty}{(bx_i)_\infty} \frac{(t x_i x_j ; q)_\infty}{(x_i x_j; q)_\infty} $$
Here $c(\lambda)$ and $r(\lambda)$ are the number of columns and rows of odd length, respectively.

\medskip

More recently some extension of the Hall-Littlewood version of the Kawanaka identity have been proved in \cite{KAWAHALL} but no proof of the more general identity has yet appeared in the litterature. The results in the second part of this thesis will eventually be published in \cite{WARNAAR}.

\medskip

We may rewrite the identity to be proved in the plethystic notation as:
$$ \sum_\lambda \Omega ( (q - \epsilon t)) B_\lambda(q,t) P_\lambda(X; q^2, t^2) = \Omega \left ( \frac{1 - \epsilon t}{1-q} X + \frac{1 - t^2}{1 - q^2} e_2(X) \right ) $$

\medskip
\noindent
By $e_2(X)$ we mean the alphabet:
$$ e_2(X) = \sum_{i < j} x_i x_j $$

\noindent
When $q=\epsilon t$ the identity reduces to the well-known generating function for the Schur functions:
$$ \sum_\lambda s_\lambda(X) = \Omega \left ( X + e_2(X) \right ) $$
which is perhaps more familiar in the form:
$$ \sum_\lambda s_\lambda(x_1, \ldots x_n) = \prod_{i = 1}^n \frac{1}{1-x_i} \prod_{1 \leq i < j \leq n} \frac{1}{1-x_i x_j} $$

\medskip
\noindent
The proof makes use of the Pieri rule and the recurrence for Macdonald polynomials, which generalize those for the Schur functions discussed in part one, followed by some combinatorial manipulations, to reduce the identity to a rational function in the variables $a_k = q^{\mu_k} t^{m - k}$. The final step is an induction on the residues at the poles.



\bigskip\bigskip
\subsection{Resultants}\label{resultants} \label{theta}

\medskip
\noindent
The {\it resultant} is defined by:
$$ R(Z : A) = \prod_{z \in Z} \prod_{a \in A} (z - a)$$
By the fundamental theorem of algebra, up to a scalar multiple, any polynomial $p(z)$ may be written in the form:
$$p(z) = R(z :A)$$
where $A$ is the alphabet of zeros. More generally, any rational function $r(z) = \frac{p(z)}{q(z)}$ may be written, up to a scalar multiple, in the form:
$$ r(z) = \frac{R(z:A)}{R(z : B)}$$
where $A$ is the alphabet of zeros $B$ is the set of poles. 

We shall need to make use of an extended version of the resultant defined by:
$$ R(X - Y : A - B) = \frac{R(X : A - B)}{R(Y : A - B)} 
= \frac{R(X : A)R(Y:B)}{R(Y : A)R(X : B)}$$
Be warned that:
$$ \Omega( X - Y ) = R( 1 : Y - X) $$

\medskip
\noindent
In the previous lemma we assumed that the alphabets $A$ and $B$ were distinct. In what follows we shall need to consider the case where one alphabet is a scalar multiple of the other.

\medskip
\noindent
Let us define:
\begin{align*}
W(X : Y)_{q,t} & =  R(X: {\color{blue}(q/t - 1)} Y) =  \prod_{\substack{x \in X \\ y \in Y}} \frac{(x - qy/t)}{(x-y)} \\
V(X : Y)_{q,t} & =  R(X: {\color{blue}(t/q - 1)} Y) = \prod_{\substack{x \in X \\ y \in Y}} \frac{(x - ty/q)}{(x-y)} 
\end{align*}
Observe that:
$$ V(X : Y)_{q,t} = W(X : Y)_{1/q, 1/t} = W(X : Y)_{t,q} $$


\medskip
\noindent
We shall also need to make use of the functions:
\begin{align*}
v(X : Y)_{q,t} & =  V(X : qY)_{q, t} = R(X: {\color{blue}(t - q)} Y) = \prod_{\substack{x \in X \\ y \in Y}} \frac{(x - ty)}{(x-qy)} \\
w(X : Y)_{q,t} & =  W(X : Y/q)_{q,t} = R(X : {\color{blue}(1/t-1/q)} Y) = \prod_{\substack{x \in X \\ y \in Y}} \frac{(x - y/t)}{(x-y/q)}
\end{align*}
Observe that:
\begin{align*}
v( X : Y)_{q,t}^{-1} & =  W( X : tY) \\
w( X : Y)_{q,t}^{-1} & =  V( X : Y/t)
\end{align*}


\bigskip\bigskip
\noindent	
Now let us define:
\begin{align*}
\Theta(X : Y)_{q,t} & =  v(X : Y)_{q,t} \,\, W(X : Y)_{q,t} \\
& =  R(X : {\color{blue}((t - q) + (q/t - 1))} Y)\\
& =  \prod_{\substack{x \in X \\ y \in Y}} \frac{(x - ty)(x - qy/t)}{(x - qy)(x - y)}
\end{align*}
as well as:
\begin{align*}
\Phi(X : Y)_{q,t} & = V(X : Y)_{q,t} \,\, w(X : Y)_{q,t} \\
& =  R(X : {\color{blue}((t/q - 1) + (1/t - 1/q))} Y)\\
& =  \prod_{\substack{x \in X \\ y \in Y}} \frac{(x - ty/q)(x - y/t)}{(x - y)(x - y/q)}
\end{align*}
Observe that:
$$ \Phi(X : Y)_{q,t} = \Theta(Y : X)_{q,t} = \Theta(X : Y)_{1/q, 1/t} = \Theta(qX : Y)_{q,t} $$

\medskip
\noindent
These rational functions will play a key role in the proof of the Kawanaka conjecture. We shall suppress the subscripts when no ambiguity can result.

\subsubsection{Residue calculations}
\begin{lemma} \label{inflim} For any $k \geq 1$ we have:
$$ \sum_{\substack{X' + X'' = X \\ |X'| = k}}\left( \Phi(X' : X'') - \Phi(X'' : X' )\right ) = 0 
$$
\begin{proof}
We wish to show that for all $1 \leq k < |A|$ we have:
\[ 
\sum_{\substack{I + J = A \\ |I| = k}} \prod_{\substack{x \in I \\ y \in J}}\frac{(x - ty/q)(x - y/t)}{(x - y)(x - y/q)}
= \sum_{\substack{I + J = A \\ |I| = k}} \prod_{\substack{x \in I \\ y \in J}}\frac{(x - qy/t)(x - ty)}{(x - y)(x - qy)}
\] 
Firstly make the substitution $ \alpha = t/q$ and $\beta = 1/t$ to rewrite this as:
\[ 
\sum_{\substack{I + J = A \\ |I| = k}} \prod_{\substack{x \in I \\ y \in J}}\frac{(x - \alpha y)(x - \beta y)}{(x - y)(x - \alpha \beta y)}
= \sum_{\substack{I + J = A \\ |I| = k}} \prod_{\substack{x \in I \\ y \in J}}\frac{(\alpha x - y)(\alpha x - y)}{(x - y)(\alpha \beta x - y)}
\] 
Write $A = x_1 + x_2 + x_3 + \cdots + x_n$. The $k=1$ case reduces to:
\[
\sum_{i=1}^n \prod_{j \neq i} \frac{(x_i - \alpha x_j)(x_i - \beta x_j)}{(x_i - x_j)(x_i - \alpha \beta x_j)}
= \sum_{i=1}^n \prod_{j \neq i} \frac{(\alpha x_i - x_j)(\beta x_i - x_j)}{(x_i - x_j)(\alpha \beta x_i - x_j)}
\]
Now consider the following contour integral:
\begin{multline*}
\frac{1}{r!(2 \pi i)^k} \oint_{C_1} \cdots \oint_{C_{r}} \frac{d w_1 \cdots d w_{r}}{w_1 \cdots w_{r}} 
\prod_{i=1}^{r} \left ( \prod_{k = 1}^n \frac{(w_i - \alpha x_k)(w_i - \beta x_k)}{(w_i - x_k)(w_i - \alpha \beta x_k)} - 1 \right ) \times \\
\prod_{i \neq j} \frac{(w_i - \alpha \beta w_j)(w_i - w_j)}{(w_i - \alpha w_j)(w_i - \beta w_j)}
\end{multline*}
where the contours of integration corresponding to the variables $w_1, w_2, \ldots w_r$ are taken such that the only poles of the integrand inside these contours are the points $x_1, \ldots x_n$. Note that the integrand has no poles at zero or at infinity.

\medskip

When $r=1$ the above integral reduces to:
\[
\frac{1}{2 \pi i} \oint_{C} \frac{dw}{w}
\left(
 \prod_{k = 1}^n \frac{(w - \alpha x_k)(w - \beta x_k)}{(w - x_k)(w - \alpha \beta x_k)}
-1\right)
\]
By Cauchy's theorem this is equal to the sum over the residues at the poles inside $C_1$ which occur precicely when $w = x_i$ for some $i$:
\[
\frac{(1 - \alpha)(1 - \beta)}{(1 - \alpha \beta)}\left ( \sum_i \prod_{j \neq i} \frac{(x_i - \alpha x_j)(x_i - \beta x_j)}{(x_i - x_j)(x_i - \alpha \beta x_j)} \right )
\]

\medskip
One can on the other hand consider the poles outside the contour: they are of the form $w = \alpha \beta x_i$ for some $i$. We obtain:
\[
\frac{(\alpha \beta - \alpha)(\alpha \beta - \beta)}{(\alpha \beta - 1)} 
\left ( \sum_i \prod_{j \neq i} \frac{(\alpha \beta x_i - \alpha x_j)(\alpha \beta x_i - \beta x_j)}{\alpha \beta (\alpha \beta x_i - x_j)(\alpha \beta x_i - \alpha \beta x_j)}\right )
\]
which simplifies to:
\[
- \frac{(1 - \alpha)(1 - \beta)}{(1 - \alpha \beta)} 
\left ( \sum_i \prod_{i \neq j} \frac{(\alpha x_i - x_j)(\beta x_i - x_j)}{(x_i - x_j)(\alpha \beta x_i - x_j)}\right )
\]
Noting that the result should be the opposite, we recover the $k=1$ version of our identity. 

\medskip

Consider now the case of general $r$. We compare once again the residues inside and outside the contours. The integration over each variable $w_i$ 
results in a sum over poles of the form $w_i=x_{a_i}$ for some $a_i$.
It is easy to check that due to the factors $w_i-w_j$, the $a_i$ must be
distinct. Furthermore, due to the symmetry of the integral by exchange of
the $w_i$, each permutation of the $a_i$ produces the same contribution
and compensates the factor $r!$. The result is a sum over subsets
$I=a_1+\cdots+a_r$ and we obtain
\begin{multline*}
\sum_{\substack{I\subset\{1,\ldots,n\}\\|I|=r}}
\frac{1}{x_{a_1}\cdots x_{a_r}}
\prod_{i=1}^r
\bigg( \frac{(x_{a_i}-\alpha x_{a_i})(x_{a_i}-\beta x_{a_i})}{(x_{a_i}-\alpha\beta x_{a_i})}
\\
\times\prod_{\substack{k=1\\k\ne a_i}}^n \frac{(x_{a_i}-\alpha x_k)(x_{a_i}-\beta x_k)}{(x_{a_i}-x_k)(x_{a_i}-\alpha\beta x_k)}
\bigg)
\prod_{i\ne j} \frac{(x_{a_i}-x_{a_j})(x_{a_i}-\alpha\beta x_{a_j})}{(x_{a_i}-\alpha x_{a_j})(x_{a_i}-\beta x_{a_j})}
\end{multline*}
which coincides with the left hand side of our identity after obvious cancellations, up to the factor
\[
((1-\alpha)(1-\beta)/(1-\alpha\beta))^r
\]

A careful analysis of the poles outside the contours shows that only the poles of the form $w_i=\alpha\beta x_{a_i}$
contribute, so that we find a similar sum:
\begin{multline*}
\sum_{\substack{I\subset\{1,\ldots,n\}\\|I|=r}}
\frac{1}{(\alpha\beta)^r x_{a_1}\cdots x_{a_r}}
\prod_{i=1}^r
\bigg( \frac{(\alpha\beta x_{a_i}-\alpha x_{a_i})(\alpha\beta x_{a_i}-\beta x_{a_i})}{(\alpha\beta x_{a_i}- x_{a_i})}
\\
\times\prod_{\substack{k=1\\k\ne a_i}}^n \frac{(\alpha\beta x_{a_i}-\alpha x_k)(\alpha\beta x_{a_i}-\beta x_k)}{(\alpha\beta x_{a_i}-x_k)(\alpha\beta x_{a_i}-\alpha\beta x_k)}
\bigg)
\prod_{i\ne j} \frac{(\alpha\beta x_{a_i}-x_{a_j})(\alpha\beta x_{a_i}-\alpha\beta x_{a_j})}{(\alpha\beta x_{a_i}-\alpha x_{a_j})(\alpha\beta x_{a_i}-\beta x_{a_j})}
\end{multline*}
which coincides with the right hand side up to the factor
\[
\frac{1}{(\alpha\beta)^r}
\left(\frac{(\alpha\beta-\alpha)(\alpha\beta-\beta)}{1-\alpha\beta}\right)^r
\]
\end{proof}
\end{lemma}

\begin{lemma}
$$ \Phi(z : X) - 1 = \frac{1-t}{1-q} \sum_{x \in X}\left( w(z : x) \Phi(X - x : x) - W(z: x) \Phi(x: X - x) \right)$$
\begin{proof}
By the $k=1$ case of the previous lemma, both sides vanish as $z$ goes to infinity, thus one need only compare the residues at the poles.
\end{proof}
\end{lemma}

\bigskip\bigskip
\begin{proposition} \label{final}
\begin{multline*}
 \sum_{s = 0}^k \frac{(q;t)_s}{(t;t)_s} \sum_{\substack{X' + X'' = X \\ |X'| = k-s}} \Big(
w(z : X'') V(z : t^{s-1} X'') W(z : t^s X') \Phi(X',X'') \\ - w(z : X') \Phi(X'',X') \Big) = 0 
\end{multline*}
\begin{proof}
When $k=1$ we have only two terms. When $s=1$ we must have $X' = \emptyset$ and $X'' = X$ thus the $s=1$ term reduces to:
$$ \frac{1-q}{1-t} \,\, (w(z :  X) V(z : X) - 1) = \frac{1-q}{1-t} \,\, (\Phi(z : X) - 1) $$
while when $s = 0$ we must have that $|X'| = 1$. Since we have:
$$ w(z : X'') V(z : X''/t) =  w(z : X'') w(z : X'')^{-1} = 1 $$
The $s=0$ term becomes:
$$ \sum_{x \in X} \left(W(z : x) \Phi(x : X - x)- w(z : x) \Phi(X - x : x)\right) $$
Thus the $k=1$ case of the identity is equivalent to the previous lemma.

\medskip
\noindent
More generally, the left hand side is a rational function in $z$ of degree zero, which vanishes in the limit as $z$ goes to infinity.
The poles are located at:
$$z \in \{a_kt^s, \,\,\,\,\, a_k/q \,\,\,:\,\,\, k = 1 \ldots m, \,\,\, s=0 \ldots k \}$$
By considering separately the terms for which $a_k \in X'$ and those for which $a_k \in X''$ one may check that the residue at the pole $z = a_k t^s$ vanishes for all $s$, while the residue at the pole $a_k / q$ is equivalent to the $k-1$ form of the identity in $m-1$ variables.
\end{proof}
\end{proposition}

\bigskip\bigskip
\subsection{Pieri formula and recurrence}\label{pierirecurrence}
\subsubsection{Arms and legs again}
For $\mu subseteq \lambda$, let $C_{\lambda / \mu}$ denote the set of boxes of $\lambda$ in columns which are longer than the corresponding columns of $\mu$ and let $R_{\lambda / \mu}$ denotes the set of boxes of $\lambda$ which are in rows longer than the corresponding rows of $\mu$. For example, if $\lambda=(8,6,5,3,2)$ and $\mu=(8,5,5,2,2)$,
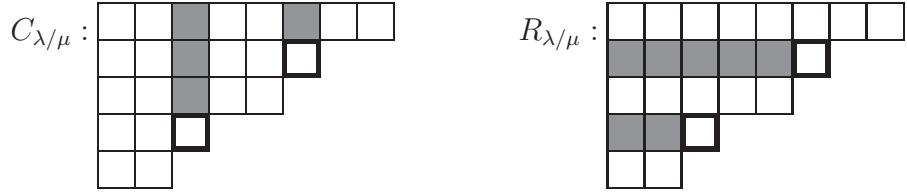
\[
C_{\lambda/\mu}:
\tableau{ \ &\ &\graycell &\ &\ &\graycell &\ &\ \\ 
\ &\ &\graycell &\ &\ &\thickcell \\ 
\ &\ &\graycell &\ &\  \\ 
\ &\ &\thickcell \\
\ &\  }
\qquad\qquad
R_{\lambda/\mu}:
\tableau{ \ &\ &\ &\ &\ &\ &\ &\ \\ 
\graycell &\graycell &\graycell &\graycell &\graycell &\thickcell \\ 
\ &\ &\ &\ &\  \\ 
\graycell &\graycell &\thickcell \\
\ &\  }
\]
Let $\tilde{C}_{\lambda / \mu}$ denote the set of boxes of $\lambda$ in columns which are the same length as the corresponding columns of $\mu$, and let $\tilde{R}_{\lambda / \mu}$ denote the set of boxes of $\lambda$ in rows which are the same length as the corresponding rows of $\mu$. With the same example,
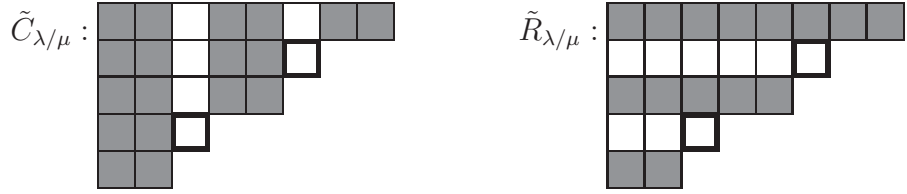
\[
\tilde C_{\lambda/\mu}:
\tableau{ \graycell&\graycell& &\graycell&\graycell& &\graycell&\graycell\\
\graycell&\graycell& &\graycell&\graycell&\thickcell \\
\graycell&\graycell& &\graycell&\graycell \\
\graycell&\graycell&\thickcell \\
\graycell&\graycell }
\qquad\qquad
\tilde R_{\lambda/\mu}:
\tableau{ \graycell&\graycell&\graycell&\graycell&\graycell&\graycell&\graycell&\graycell\\
 & & & & &\thickcell \\
\graycell&\graycell&\graycell&\graycell&\graycell \\
 & &\thickcell \\
\graycell&\graycell }
\]

\medskip
Now let us define:
$$ C_{\lambda / \mu}(q,t) = \sum_{s \in C_{\lambda / \mu}} q^{a_\lambda(s)} t^{l_\lambda(s)} - \sum_{s\in \mu\cap C_{\lambda/\mu}}q^{a_\mu(s)} t^{l_\mu(s)} $$
$$ R_{\lambda / \mu}(q,t) = \sum_{s \in R_{\lambda / \mu}} q^{a_\lambda(s)} t^{l_\lambda(s)} - \sum_{s\in \mu\cap R_{\lambda/\mu}}q^{a_\mu(s)} t^{l_\mu(s)} $$
$$ \tilde{C}_{\lambda / \mu}(q,t) = \sum_{s \in \tilde C_{\lambda / \mu}} (q^{a_\lambda(s)} t^{l_\lambda(s)} - q^{a_\mu(s)} t^{l_\mu(s)}) $$
$$ \tilde{R}_{\lambda / \mu}(q,t) = \sum_{s \in \tilde R_{\lambda / \mu}} (q^{a_\lambda(s)} t^{l_\lambda(s)} - q^{a_\mu(s)} t^{l_\mu(s)}) $$

\medskip
\noindent
It is clear that we have:
$$ C_{\lambda' / \mu'}(t,q) = R_{\lambda / \mu}(q,t) $$
$$ \tilde{C}_{\lambda' / \mu'}(t,q) = \tilde{R}_{\lambda / \mu}(q,t) $$
as well as that:
\begin{align*}
B_\lambda(q,t) - B_\mu(q,t)  
& = C_{\lambda / \mu}(q,t) + \tilde{C}_{\lambda / \mu}(q,t)  \\
& = R_{\lambda / \mu}(q,t) + \tilde{R}_{\lambda / \mu}(q,t) 
\end{align*}


\bigskip\bigskip
\subsubsection{Pieri formula}
The Macdonald polynomials are, in fact, the eigenfunctions of the more general family of operators \cite{BLACK} (page 315):
$$ \frac{1}{\Delta_n(X)} \sum_{\substack{I \subseteq [n] \\ |I| = r}} \Delta_n(X + (1-t)X_I) f(X + (1-q)X_I)$$
where $X_I$ is the alphabet:
$$ X_I = \sum_{i \in I} x_i$$
with eigenvalues:
$$e_r \left (\sum_{i=1}^n q^{\lambda_i} t^{n-i} \right )$$

\medskip
Using this fact, Macdonald was also able to show that his polynomials satisfy the following generalization of the Pieri formulae \cite{BLACK} (page 340):

\begin{align*}
P_\mu(X)  \Omega_z \left(X \frac{1-t}{1-q}\right) & =  \sum_{\lambda \in U(\mu)} \varphi_{\lambda / \mu}(q,t) P_\lambda(q,t) z^{|\lambda| - |\mu|} \\
Q_\mu(X) \tilde{\Omega}_z(X) & = \sum_{\lambda \in \tilde{U}(\mu)} \varphi'_{\lambda / \mu}(q,t) Q_\lambda(X) z^{|\lambda| - |\mu|} \\
\end{align*}
where:
\begin{align*}
\varphi_{\lambda / \mu}(q,t) & = \Omega ( (q - t) C_{\lambda / \mu}(q,t)) \\
\varphi'_{\lambda / \mu}(q,t) & = \Omega ( (t - q) R_{\lambda / \mu}(q,t) ) \\
\end{align*}
As before $U(\mu)$ denotes the set of partitions obtained from $\mu$ by adding a horizontal strip while $\tilde{U}(\mu)$ denotes the set of partitions obtained from $\mu$ by adding a vertical strip. 
Note that the latter may be obtained from the former by applying the operator $\omega_{q,t}$ to both sides and then formally replacing $\lambda$ by $\lambda'$ and $\mu$ by $\mu'$. 
The two additional Pieri formulae:
\begin{align*}
Q_\mu(X) \Omega_z \left(X \frac{1-t}{1-q}\right) & =  \sum_{\lambda \in U(\mu)} \psi_{\lambda / \mu}(q,t) Q_\lambda(X)z^{|\lambda| - |\mu|} \\
P_\mu(X) \tilde{\Omega}_z(X) & = \sum_{\lambda \in \tilde{U}(\mu)}\psi'_{\lambda / \mu}(q,t) P_\lambda(X) z^{|\lambda| - |\mu|}
\end{align*}
where:
\begin{align*}
\psi_{\lambda / \mu}(q,t) & =  \Omega((t - q) \tilde{C}_{\lambda / \mu}(q,t)) \\
\psi'_{\lambda / \mu}(q,t) & =  \Omega((q - t) \tilde{R}_{\lambda / \mu}(q,t)) \\
\end{align*}
are obtained from by multiplying, or dividing by $\Omega((q-t)(B_\lambda(q,t) - B_\mu(q,t)))$ as appropriate.

\medskip
\noindent
\subsubsection{Recurrence}
Dual to the Pieri formulae, we have the following recurrences \cite{BLACK} (pages 346, 348):
$$ P_\lambda(X + z) = \sum_{\mu \in D(\lambda)}\psi_{\lambda / \mu}(q,t)P_\mu(X) z^{|\lambda| - |\mu|}$$
$$ Q_\lambda(X - z) = \sum_{\mu \in \tilde{D}(\lambda)}\psi'_{\lambda / \mu}(q,t)Q_\mu(X) (\epsilon z)^{|\lambda| - |\mu|}$$
$$ Q_\lambda(X + z) = \sum_{\mu \in D(\lambda)} \varphi_{\lambda / \mu}(q,t)Q_\mu(X) z^{|\lambda| - |\mu|}$$
$$ P_\lambda(X - z) = \sum_{\mu \in \tilde{D}(\lambda)} \varphi'_{\lambda / \mu}(q,t) P_\mu(X)(\epsilon z)^{|\lambda| - |\mu|}$$
Here $D(\lambda)$ denotes the set of partitions obtained from $\lambda$ by removing a horizontal strip, while $\tilde{D}(\lambda)$ denotes the set of partitions obtained from $\lambda$ by removing a vertical strip.

\smallskip
Note that it is also possible to take any one of the recurrences or the Pieri formulae, together with the definition of the Macdonald polynomial of a row, or a column accordingly  as the {\it definition} of the Macdonald polynomials. 

\bigskip\bigskip
\subsection{The Proof}\label{proof}
\subsubsection{The Schur case}
To aid the reader in not getting lost in a thicket of $q$'s and $t$'s we first give a proof of the classical Schur identity, which mirrors, in simplified form, the main steps of the proof of the more complicated Kawanaka identity. 
\begin{proposition} \label{schur}
$$ \sum_\lambda s_\lambda(x_1, \ldots, x_{n+1})= \prod_{i=1}^{n} \frac{1}{(1-x_i)} \prod_{1 \leq i<j \leq n} \frac{1}{(1-x_i x_j)} $$
\begin{proof}
The proof is by induction on the number of variables. The first step makes use of the Pieri formula and recurrence for Schur functions. The second step is a bijection between two families of partitions. 

When $n=1$ we get the expansion of the geometric series:
$$ \sum_{k \geq 0} x_1^k = \frac{1}{1-x_1} $$
By the recurrence for the Schur functions we have:
\begin{align*}
\text{LHS}(n+1) & = \sum_\lambda s_\lambda(x_1, \ldots, x_{n+1}) \\
& =  \sum_\lambda  \sum_{\mu \in D(\lambda)} 
s_\mu(x_1, \ldots, x_{n}) x_{n+1}^{|\lambda| - |\mu|} \\
& = \sum_\mu s_\mu(x_1, \ldots, x_{n}) 
\sum_{\lambda \in U(\mu)} x_{n+1}^{|\lambda| - |\mu|} 
\end{align*}
By the induction hypothesis:
\begin{align*}
\text{RHS}(n+1) & = 	\prod_{i=1}^{n+1} \frac{1}{1-x_i} \prod_{1 \leq i<j \leq n+1} \frac{1}{1-x_i x_j} \\
& =  \frac{1}{1-x_{n+1}} 
 \prod_{i=1}^{n} \frac{1}{1-x_i} \prod_{1 \leq i<j \leq n} \frac{1}{1-x_i x_j} \prod_{i=1}^{n} \frac{1}{1-x_i x_{n+1}} \\
& =  \frac{1}{1-x_{n+1}} 
 \text{RHS}(n) \prod_{i=1}^{n} \frac{1}{1-x_i x_{n+1}} \\
& =  \frac{1}{1-x_{n+1}} 
\text{LHS}(n)
\prod_{i=1}^{n} \frac{1}{1-x_i x_{n+1} } \\
& =  \frac{1}{1-x_{n+1}}  
\sum_\mu s_\mu(x_1, \ldots, x_{n}) \prod_{i=1}^{n}\frac{1}{1-x_i x_{n+1} }
\end{align*}
But by the Pieri formula for the Schur functions:	
\begin{align*}
\sum_\mu s_\mu(x_1, \ldots, x_{n}) \prod_{i=1}^{n}\frac{1}{1-x_i x_{n+1} } & =  
\sum_\mu s_\mu(x_1, \ldots, x_{n}) 
\sum_s h_s(x_1, \ldots, x_{n}) x_{n+1}^s \\
& =  \sum_\mu  \sum_{\lambda \in U(\mu)}
s_\lambda(x_1, \ldots, x_{n}) x_{n+1}^{|\lambda| - |\mu|} \\
& =  \sum_\lambda s_\lambda(x_1, \ldots, x_{n}) 
\sum_{\mu \in D(\lambda)} x_{n+1}^{|\lambda| - |\mu|} 
\end{align*}

\bigskip\bigskip
\noindent
By considering the coefficient of $s_\mu(X)$ on both sides, it suffices to demonstrate that:
\begin{align}\label{part1}
\sum_{\lambda \in U(\mu)} x_{n+1}^{|\lambda| - |\mu|}
& =  \frac{1}{1-x_{n+1}} \sum_{\gamma \in D(\mu)} x_{n+1}^{|\mu| - |\gamma|}
\end{align}
Conjugating all the partitions involved, this is equivalent to:
\begin{align}\label{part2}
\sum_{\lambda \in \tilde{U}(\mu)} x_{n+1}^{|\lambda| - |\mu|}
& =  \frac{1}{1-x_{n+1}} \sum_{\gamma \in \tilde{D}(\mu)} x_{n+1}^{|\mu| - |\gamma|}
\end{align}

\medskip
Suppose now that $\mu$ is a partition with distinct $m$ distinct nonzero parts.
For each $\alpha \subseteq [m]$ such that $|\alpha| = k$ let $\mu_-(\alpha)$ denote the partition obtained from $\mu$ by removing a box from the end of each row indexed by an element of $\alpha$. We have:
$$\tilde{D}_k(\mu) = \{ \mu_-(\alpha), \alpha \subseteq [m], |\alpha| = k \}$$

Similarly, let $\mu_+(\alpha)$ denote the partition obtained from $\mu$ by adding a box to the end of each row indexed by $\alpha$. Let:
$$S\tilde{U}_k(\mu) = \{ \mu_+(\alpha), \alpha \subseteq [m], |\alpha| = k \}$$

There is a natural bijection between $\tilde{D}_k(\mu)$ and $S\tilde{U}_k(\mu)$ which send $\mu_-(\alpha)$ to $\mu_+(\alpha)$. 

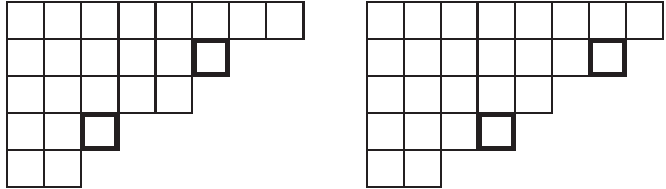
\[
\tableau{ \ &\ &\ &\ &\ &\ &\ &\ \\ 
\ &\ &\ &\ &\ &\thickcell \\ 
\ &\ &\ &\ &\  \\ 
\ &\ &\thickcell \\
\ &\  }
\,\,\,\,\,\,\,\,\,\,\,\,
\tableau{ \ &\ &\ &\ &\ &\ &\ &\ \\ 
\ &\ &\ &\ &\ &\ &\thickcell \\ 
\ &\ &\ &\ &\  \\ 
\ &\ &\ &\thickcell \\
\ &\  }
\]

\medskip
If $\mu$ contains repeated parts, then it is no longer true that $\mu_+(\alpha)$ and $\mu_-(\alpha)$ will be valid partitions for all $\alpha$. 
For example, if $\mu_i = \mu_{i+1}$ and $i \in \alpha$ but $(i+1) \not\in \alpha$ then $\mu_-(\alpha)$ will not be a valid partition. 

Conversely, if $i \not \in \alpha$ while $(i+1) \in \alpha$ then $\mu_+(\alpha)$ will not be a valid partition. 
There is still, however a natural bijection between $\tilde{D}_k(\mu)$ and $S\tilde{U}_k(\mu)$, though it is slightly more cumbersome to describe.

\medskip
More generally for each $\alpha \subseteq [m]$ such that $|\alpha| = k$ and for each $p \geq 0$ let $\mu_+(\alpha,p)$ denote the partition obtained from $\mu$ by adding a box to the end of each row indexed by an element of $\alpha$, and then adding $p$ new rows each of length one. Let:
$$S\tilde{U}_k(\mu,p) = \{ \mu_+(\alpha,p), \alpha \subseteq [m], |\alpha| = k \}$$

Upon consideration of the coefficient of $x_n^k$ in equation \eqref{part2}, to complete the proof we must show that:
\begin{equation} \label{bijection}
|\tilde{U}_k(\mu)| = \sum_s |\tilde{D}_{k-s}(\mu)| 
\end{equation}
This follows immediately from the fact that:
$$ \tilde{U}_k(\mu) = \bigcup_{p=0}^k S\tilde{U}_{k-p}(\mu,p) $$ 
and that for each $p$ we have:
$$ |S\tilde{U}_k(\mu,p)| = |S\tilde{U}_k(\mu,0)| = |\tilde{D}_k(\mu)| $$

\end{proof}
\end{proposition}

\subsubsection{Step one}
\noindent
For the Kawanaka identity, let us define:
\begin{align*}
\text{LHS}(n) & =  \sum_{ \lambda} \Omega \left ( {\color{blue}(q - \epsilon t) B_\lambda(q,t)}\right ) P_\lambda(X_n;q^2,t^2) \\
\text{RHS}(n) & =  \Omega \left ( {\color{blue}\frac{1-\epsilon t}{1-q}} X_n + {\color{blue}\frac{1-t^2}{1-q^2}} e_2(X_n) \right ) 
\end{align*}
	


\medskip
\noindent
The case $n=1$ is a special case of the $q$-binomial formula:
\begin{align*}
\Omega \left ( {\color{blue}\frac{1-\epsilon t}{1 - q}}x_1\right ) 
&= \sum_n \Omega \left ( {\color{blue}\frac{(q-\epsilon t)(1-q^n)}{(1-q)}}\right ) x_1^n \\
&= \sum_n \Omega \left ( {\color{blue}(q-\epsilon t) B_{(n)}(q, \epsilon t)} \right ) x_1^n
\end{align*}


\medskip
\noindent
We have now, by the recurrence for the Macdonald $P$-functions:
\begin{align*}
&\text{LHS}(n+1) 
\\&= \sum_{ \lambda} \Omega \left ( {\color{blue}(q - \epsilon t) B_\lambda(q,t)}\right ) P_\lambda(X_n + x_{n+1};q^2,t^2) \\
&= \sum_{ \lambda} \Omega \left ( {\color{blue}(q - \epsilon t) B_\lambda(q,t)}\right ) \!\!\sum_{\mu \in D(\lambda) }\!\!\!\Omega({\color{blue}(t^2 - q^2) \tilde{C}_{\lambda / \mu}(q^2,t^2)}) P_\mu(X_n;q^2, t^2) x_{n+1}^{|\lambda| - |\mu|} \\
&= \sum_\mu P_\mu(X_n,q^2,t^2) \sum_k x_{n+1}^k\!\!\! \sum_{\lambda \in U_k(\mu)}\!\!\!\! \Omega \left ( {\color{blue}(q - \epsilon t) B_\lambda(q,t)}\right ) \Omega({\color{blue}(t^2 - q^2) \tilde{C}_{\lambda / \mu}(q^2,t^2)})
\end{align*}

\medskip
\noindent
while, by the induction hypothesis:
\begin{align*}
\text{RHS}(n+1) &= \Omega \left ( {\color{blue}\frac{1-\epsilon t}{1-q}} x_{n+1}\right ) \text{RHS}(n)  \Omega \left ({\color{blue}\frac{1-t^2}{1-q^2}} X_n x_{n+1} \right ) \\
&= \Omega \left ( {\color{blue}\frac{1-\epsilon t}{1-q}} x_{n+1}\right ) \text{LHS}(n) \Omega \left ({\color{blue}\frac{1-t^2}{1-q^2}} X_n x_{n+1} \right )  \end{align*}

\medskip
\noindent
But, by the Pieri formula for the Macdonald $P$ functions, we have:
\begin{align*}
& \text{LHS}(n) \Omega \left ({\color{blue}\frac{1-t^2}{1-q^2}} X_n x_{n+1} \right ) \\
&= \sum_{ \gamma} \Omega \left ( {\color{blue}(q - \epsilon t) B_\gamma(q,t)}\right ) P_\gamma(X_n;q^2,t^2)\,\Omega \left ({\color{blue}\frac{1-t^2}{1-q^2}} X_n x_{n+1} \right ) \\
&= \sum_{ \gamma}  \Omega \left ( {\color{blue}(q - \epsilon t) B_\gamma(q,t)}\right ) \!\!\sum_{\mu \in U(\gamma)} \!\!\Omega ( {\color{blue}(q^2 - t^2) C_{\mu / \gamma}(q^2,t^2)}) P_\mu(X_n;q^2,t^2) x_{n+1}^{|\mu| - |\gamma|}  \\
&= \sum_\mu P_\mu(X_n;q^2,t^2) \sum_k x_{n+1}^k \!\!\!\sum_{\gamma \in D_k(\mu)}\!\!\!\! \Omega \left ( {\color{blue}(q - \epsilon t) B_\gamma(q,t)}\right ) \Omega ( {\color{blue}(q^2 - t^2) C_{\mu / \gamma}(q^2,t^2)} )
\end{align*}

\bigskip\bigskip
\noindent
By considering the coefficient of $P_\mu(X;q^2,t^2)$ on both sides, and then dividing by $\Omega({\color{blue}(q-\epsilon t)B_\mu(q,t)})$, we need to show that, for any $\mu$ we have:
\begin{align*}
&\sum_k x_{n+1}^k \sum_{\lambda \in U_k(\mu)} \Omega \left ( {\color{blue}(q - \epsilon t) (B_\lambda(q,t) - B_\mu(q,t))}\right ) \Omega({\color{blue}(t^2 - q^2) \tilde{C}_{\lambda / \mu}(q^2,t^2)}) \\
&= \Omega \left ( {\color{blue}\frac{1-\epsilon t}{1-q}} x_{n+1} \right ) \times
\\
&\sum_k x_{n+1}^k \sum_{\gamma \in D_k(\mu)} \Omega \left ( {\color{blue}(q - \epsilon t) (B_\gamma(q,t) - B_\mu(q,t))}\right ) \Omega ( {\color{blue}(q^2 - t^2) C_{\mu / \gamma}(q^2,t^2)})
\end{align*}
Observe that when $q = \epsilon t$ this reduces to equation \eqref{part1}. 

\medskip
\noindent
By conjugating all the partitions involved, and interchanging the roles of $q$ and $t$ we obtain the dual form:

\begin{align*}
&\sum_k x_{n+1}^k \sum_{\lambda \in \tilde{U}_r(\mu)} \Omega({\color{blue}(t-\epsilon q) (B_\lambda(q,t) - B_\mu(q,t))}) \Omega({\color{blue}(q^2 - t^2) \tilde{R}_{\lambda / \mu}(q^2, t^2)}) \\
&= \Omega \left ( {\color{blue}\frac{1-\epsilon q}{1-t}} x_{n+1} \right ) \times \\
&\sum_k x_{n+1}^k  \sum_{\gamma \in \tilde{D}_z(\mu)} \Omega({\color{blue}(t-\epsilon q)(B_\gamma(q,t) - B_\mu(q,t))}) \Omega({\color{blue}(t^2 - q^2) R_{\lambda / \mu}(q^2,t^2)})
\end{align*}


\medskip
Here we have made use of the fact that ${\color{blue}B_{\lambda'}(t,q) = B_\lambda(q,t)}$ as well as ${\color{blue}\tilde{C}_{\lambda' / \mu'}(t,q) = \tilde{R}_{\lambda / \mu}(q,t)}$ and ${\color{blue}C_{\lambda' / \mu'}(t,q) = R_{\lambda / \mu}(q,t)}$. 

\medskip
\noindent
Note that when $q = \epsilon t$ this reduces to equation \eqref{part2}.

\medskip
\noindent
Finally, making use of the $q$-binomial formula once again, on consideration of the coefficient of $x_{n+1}^r$ on both sides, we must show that for any $\mu$ and any $r$ we have:
\begin{equation} \label{ident}\begin{split}
 &\sum_{\lambda \in \tilde{U}_r(\mu)} \Omega \left ( {\color{blue}(t - \epsilon q) (B_\lambda(q,t)  - B_\mu(q,t))} \right ) \Omega ( {\color{blue}(q^2 - t^2) R_{\lambda / \mu}(q^2,t^2)}) \\
&\quad = \sum_s \Omega({\color{blue}(t-\epsilon q)B_{(1^{r-s})}(\epsilon q, t)}) \times \\
&\  \sum_{\gamma \in \tilde{D}_s(\mu)} \Omega \left ( {\color{blue}(t - \epsilon q) (B_\gamma(q,t) - B_\mu(q,t))} \right ) \Omega({\color{blue}(t^2 - q^2) \tilde{R}_{\lambda / \mu}(q^2,t^2)})
\end{split}\end{equation}

\medskip
\noindent
Let us define:
\begin{align*}
H_\lambda(s) &= \Omega({\color{blue}(q-\epsilon t) q^{a_\lambda(s)} t^{l_\lambda(s)}}) 
= \frac{1+q^{a_\lambda(s)} t^{l_\lambda(s)+1}}{1 - q^{a_\lambda(s)+1} t^{l_\lambda(s)}}\\
\tilde{H}_\lambda(s) &= \Omega({\color{blue}(t-\epsilon q) q^{a_\lambda(s)} t^{l_\lambda(s)}}) 
= \frac{1+q^{a_\lambda(s)+1} t^{l_\lambda(s)}}{1 - q^{a_\lambda(s)} t^{l_\lambda(s)+1}}\\
G_\lambda(s) &= \Omega({\color{blue}(t - q) q^{2a_\lambda(s)} t^{2l_\lambda(s)}}) 
= \frac{1-q^{2(a_\lambda(s)+1)} t^{2l_\lambda(s)}}{1 - q^{2a_\lambda(s)} t^{2(l_\lambda(s)+1)}} \\
\end{align*}
We have, by the difference of perfect squares $\Omega({\color{blue}a^2}) = \Omega({\color{blue}a + \epsilon a})$ that:

\begin{align*}
G_\lambda(s) H_\lambda(s) &= \Omega({\color{blue}(t^2-q^2) q^{2a_\lambda(s)} t^{2l_\lambda(s)} + (q-\epsilon t) q^{a_\lambda(s)} t^{l_\lambda(s)}}) \\
&= \Omega({\color{blue}(t-\epsilon q) q^{a_\lambda(s)} t^{l_\lambda(s)}}) \\
&= \tilde{H}_\lambda(s)
\end{align*}

\bigskip\bigskip
\noindent
Now we may simplify equation \eqref{ident} by writing:
\begin{align*}
L(\lambda, \mu) &= \Omega \left ( {\color{blue}(t - \epsilon q) (B_\lambda(q,t) - B_\mu(q,t)} \right ) \Omega({\color{blue}(q^2 - t^2) \tilde{R}_{\lambda / \mu}(q^2,t^2)})  \\
&= \Omega \left ( {\color{blue}(q-\epsilon t)\tilde{R}_{\lambda / \mu}(q,t)  + (t - \epsilon q)R_{\lambda / \mu}(q,t)} \right )  \\
&=  \prod_{s \in \tilde{R}_{\lambda / \mu}} \frac{H_\lambda(s)}{H_\mu(s)} \prod_{s \in R_{\lambda / \mu}} \frac{\tilde{H}_\lambda(s)}{\tilde{H}_\mu(s)} \\
&=  \prod_{s \in \tilde{R}_{\lambda / \mu} \cap C_{\lambda / \mu}} \frac{H_\lambda(s)}{H_\mu(s)} \prod_{s \in R_{\lambda / \mu}} \frac{\tilde{H}_\lambda(s)}{\tilde{H}_\mu(s)} 
\end{align*}

\medskip
\noindent
and:
\begin{align*}
R(\mu, \gamma) &= \Omega \left ( {\color{blue}(t - \epsilon q) (B_\gamma(q,t) - B_\mu(q,t))} \right ) \Omega ( {\color{blue}(t^2 - q^2) \tilde{R}_{\mu / \gamma}(q^2,t^2)}) \\
&=\Omega \left ( {\color{blue}-(q - \epsilon t) R_{\mu / \gamma}(q,t) - (t - \epsilon q) \tilde{R}_{\mu / \gamma}(q,t)}  \right ) \\
&= \prod_{s \in R_{\mu / \gamma}} \frac{H_\gamma(s)}{H_\mu(s)} \prod_{s \in \tilde{R}_{\mu / \gamma}} \frac{\tilde{H}_\gamma(s)}{\tilde{H}_\mu(s)} \\
&= \prod_{s \in R_{\mu / \gamma}} \frac{H_\gamma(s)}{H_\mu(s)} \prod_{s \in \tilde{R}_{\mu / \gamma} \cap C_{\mu / \gamma}} \frac{\tilde{H}_\gamma(s)}{\tilde{H}_\mu(s)} 
\end{align*}

\medskip
\noindent
To complete the proof we must show that for any partition $\mu$ and for any positive integer $k$ we have:
\begin{align} \label{toprove}
\sum_{\lambda \in \tilde{U}_k(\mu)} L(\lambda, \mu) = \sum_s \frac{(-q;t)_s}{(t;t)_s} \sum_{\gamma \in \tilde{D}_{k-s}(\mu)} R(\mu, \gamma) 
\end{align}
Note that when $q = \epsilon t$ we have $L(\lambda, \mu) = R(\mu, \gamma) = 1$ and the above reduces to equation \eqref{bijection}.


\subsubsection{Step two}
Suppose that $\mu$ is a partition with $m$ distinct nonzero parts. For each $k \in [m] = \{1,2, \ldots m \}$ let $a_k = q^{\mu_k} t^{m-k}$. For any $K \subseteq [m]$ let $A_K$ denote the alphabet:
$$ A_K = \sum_{k \in K} a_k $$
The next step is to obtain several related combinatorial descriptions for the rational function first introduced in section \ref{theta}: 
$$\Phi(A_I,A_J)_{\epsilon q, t} = \prod_{\substack{i \in I \\ j \in J}} \frac{(a_i + t/q \,\, a_j)(a_i -  a_j/t)}{(a_i - a_j)(a_i + a_j /q)}$$ 
Note that $q$ has been replaced by $\epsilon q$, that is every occurence of $q$ has been replaced by its formal (as opposed to plethystic) negative.

\bigskip\bigskip
Let $\alpha$ and $\tilde{\alpha}$ be such that $\alpha \cup \tilde{\alpha} = [m]$. Let $\gamma = \mu_-(\alpha)$ denote the partition obtained from $\mu$ by removing a box from the end of each row indexed by $\alpha$ and let $\tilde{\gamma} = \mu_-(\tilde{\alpha})$ denote the partition obtained from $\mu$ by removing a box from the end of each row indexed by $\tilde{\alpha}$. 

\medskip
\noindent
If $s \in \tilde{R}_{\mu / \gamma} \cap {C}_{\mu / \gamma}$ then $s$ must be of the form $s = (i,\mu_j)$ for some $j \in \alpha$ and $i \not\in \alpha$ with $i < j$. In this case we have that $a_{\gamma}(s) = a_{\mu}(s) = \mu_i - \mu_j$,
while
$ l_\gamma(s) +1 = l_{\mu}(s) = j - i $.

\medskip
\noindent
For example: 
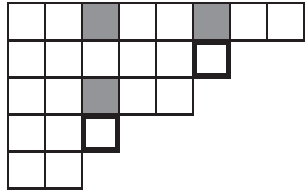
\[
\tableau{ \ &\ &\graycell &\ &\ &\graycell &\ &\ \\ 
\ &\ &\ &\ &\ &\thickcell \\ 
\ &\ &\graycell &\ &\  \\ 
\ &\ &\thickcell \\
\ &\  }
\]
$\mu = (8,6,5,3,2)$ and $\alpha = \{2,4 \}$, $\tilde{R}_{\mu / \gamma} \cap {C}_{\mu / \gamma} = \{ (1, \mu_4), (1, \mu_2), (3, \mu_4) \} = \{ (1, 3), (1, 6), (3, 3) \}$.

\medskip
\noindent
In this case, one may check that:
$$ \tilde{H}_\mu(i, \mu_j) = \frac{1 + q \, a_i a_j^{-1}}{1 - t \, a_i a_j^{-1}} = v(a_j : a_i)_{\epsilon q, t}^{-1} $$
$$ \tilde{H}_{\mu_-(\alpha)}(i, \mu_j) =  \frac{1 + q /t \,\, a_i a_j^{-1}}{1 - a_i a_j^{-1}} = W(a_j : a_i)_{\epsilon q, t} $$
so that:
\begin{equation} \label{one}
\frac{\tilde{H}_{\mu_-(\alpha)}(i,\mu_j)}{\tilde{H}_\mu(i,\mu_j)} = \Theta(a_j,a_i)_{\epsilon q, t} = \Phi(a_i,a_j)_{\epsilon q, t} 
\end{equation}

\medskip
\noindent
Similarly if $s \in \tilde{R}_{\mu / \tilde{\gamma}} \cap {C}_{\mu / \tilde{\gamma}}$ then $s$ must be of the form $s = (j,\mu_i)$ for some $j \not\in \tilde{\alpha}$ and $i \in \tilde{\alpha}$ with $j < i$.
now $a_{\tilde{\gamma}}(s) = a_{\mu}(s) = \mu_j - \mu_i$,
while
$ l_{\tilde{\gamma}}(s) +1 = l_{\mu}(s) = i - j $.

\medskip
\noindent
For example:
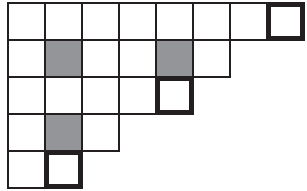
\[
\tableau{ \ &\ &\ &\ &\ &\ &\ &\thickcell \\ 
\ &\graycell &\ &\ &\graycell &\ \\ 
\ &\ &\ &\ &\thickcell  \\ 
\ &\graycell &\ \\
\ &\thickcell  }
\]
$\mu = (8,6,5,3,2)$ and $\tilde{\alpha} = \{1,3,5 \}$, $\tilde{R}_{\mu / \tilde{\gamma}} \cap {C}_{\mu / \tilde{\gamma}} = \{ (2, \mu_5), (2, \mu_3), (4, \mu_5) \} = \{ (2, 2), (2, 5), (4, 2) \}$.

\medskip
\noindent
In this case we have:
\begin{align*}
\tilde{H}_\mu(j, \mu_i) &= v(a_i : a_j)_{\epsilon q, t}^{-1} \\
\tilde{H}_{\mu_-(\tilde{\alpha})}(j, \mu_i) &= W(a_i : a_j)_{\epsilon q,  t}
\end{align*}
so that:
\begin{equation}\label{two}
\frac{\tilde{H}_{\mu_-(\tilde{\alpha})}(j, \mu_i)}{\tilde{H}_\mu(j, \mu_i)} = \Theta(a_i, a_j)_{\epsilon q, t} = \Phi(a_j, a_i)_{\epsilon q, t}
\end{equation}

\bigskip\bigskip
\noindent
Now let $\lambda = \mu_+(\alpha)$ and let $\tilde{\lambda} = \mu_+(\tilde{\alpha})$. 

\medskip
\noindent
If $s \in \tilde{R}_{\lambda / \mu} \cap {C}_{\lambda / \mu}$ then we must have $s = (i,\mu_j+1)$ for some $j \in \alpha$ and $i \not\in \alpha$ with $i < j$. This time we have that $a_{\mu}(s) = a_{\lambda}(s) = \mu_i - \mu_j - 1$,
while
$ l_\mu(s) +1 = l_{\lambda}(s) = j - i $.

\medskip
\noindent
For example:
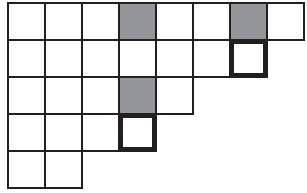
\[
\tableau{ \ &\ &\ &\graycell &\ &\ &\graycell &\  \\ 
\ &\ &\ &\ &\ &\ &\thickcell \\ 
\ &\ &\ &\graycell &\  \\ 
\ &\ &\ &\thickcell \\
\ &\  }
\]
If $\mu = (8,6,5,3,2)$ and $\alpha = \{2,4 \}$ then $\tilde{R}_{\lambda / \mu} \cap {C}_{\lambda / \mu}	 = \{ (1, \mu_4+1), (1, \mu_2+1), (3, \mu_4+1) \} = \{ (1, 4), (1, 7), (3, 4) \}$.

\medskip
\noindent
Now:
$$ H_{\mu_+(\alpha)}(i, \mu_j+1) = \frac{1 + t/q \,\, a_i a_j^{-1}}{1 - a_i a_j^{-1}} = V(a_j, a_i)_{\epsilon q, t} $$
$$ H_\mu(i, \mu_j + 1) = \frac{1 + 1/q \,\, a_i a_j^{-1}}{1 - 1/t \,\, a_i a_j^{-1}} = w(a_j : a_i)_{ \epsilon q, t}^{-1} $$
so that:
\begin{equation} \label{three}
 \frac{H_{\mu_+(\alpha)}(i, \mu_j+1)}{H_\mu(i, \mu_j + 1)} = \Phi(a_j, a_i)_{\epsilon q, t} 
\end{equation}

\medskip
\noindent
Similarly, if $s \in \tilde{R}_{\tilde{\lambda} / \mu} \cap {C}_{\tilde{\lambda} / \mu}$ then we must have $s = (j,\mu_i+1)$ for some $j \not\in \tilde{\alpha}$ and $i \in \tilde{\alpha}$ with $j < i$. This time $a_{\mu}(s) = a_{\lambda}(s) = \mu_j - \mu_i - 1$,
while
$ l_\mu(s) +1 = l_{\lambda}(s) = i - j $

\medskip
\noindent
For example:
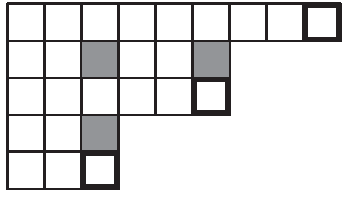
\[
\tableau{ \ &\ &\ &\ &\ &\ &\ &\ &\thickcell  \\ 
\ &\ &\graycell &\ &\ &\graycell \\ 
\ &\ &\ &\ &\ &\thickcell \\ 
\ &\ &\graycell \\
\ &\ &\thickcell }
\]
If $\mu = (8,6,5,3,2)$ and $\tilde{\alpha} = \{1,3,5 \}$ then $\tilde{R}_{\tilde{\lambda} / \mu} \cap {C}_{\tilde{\lambda} / \mu} = \{ (2, \mu_5+1), (2, \mu_3+1), (4, \mu_5+1) \} = \{ (2, 3), (2, 6), (4, 3) \}$

\medskip
\noindent
We have:
\begin{align*}
 H_{\mu_+(\tilde{\alpha})}(j, \mu_i+1) &=  V(a_i, a_j)_{\epsilon q, t} \\
H_\mu(j, \mu_i + 1) &= w(a_i, a_j)_{\epsilon q,  t}^{-1}
\end{align*}
so that:
\begin{equation} \label{four}
\frac{H_{\mu_+(\tilde{\alpha})}(j, \mu_i+1)}{H_\mu(j, \mu_i + 1)} = \Phi(a_i, a_j)_{\epsilon q, t} 
\end{equation}	

\bigskip\bigskip
\noindent
Note that in equations \eqref{one} and \eqref{three} we have the restriction that $i<j$ while in equations \eqref{two} and \eqref{four} we have the restriction that $j < i$.

\medskip
\noindent
We may combine all these observations into the following:

\begin{proposition} 
We have:
\begin{align*}
\Phi(A_I,A_J)_{\epsilon q, t} &= \prod_{\substack{i \not\in \alpha \\ j \in \alpha \\ i < j}} \frac{\tilde{H}_{\mu_-(\alpha)}(i,\mu_j)}{\tilde{H}_\mu(i,\mu_j)}
\prod_{\substack{i \not\in \alpha \\ j \in \alpha \\ j < i}} \frac{H_{\mu_+(\tilde{\alpha})}(j,\mu_i+1)}{H_\mu(j,\mu_i+1)} \\
\end{align*}
and:
\begin{align*}
\Phi(A_J,A_I)_{\epsilon q, t} &= \prod_{\substack{i \not\in \alpha \\ j \in \alpha \\ i < j}} \frac{H_{\mu_+(\alpha)}(i,\mu_j+1)}{H_\mu(i,\mu_j+1)} 
\prod_{\substack{i \not\in \alpha \\ j \in \alpha \\ j < i}} \frac{\tilde{H}_{\mu_-(\tilde{\alpha})}(j,\mu_i)}{\tilde{H}_\mu(j,\mu_i)} \\
\end{align*}

\begin{proof}
The first expression is obtained by combining equations \eqref{one} and \eqref{four}. The second expression is obtained by combining equations \eqref{two} and \eqref{three}.
\end{proof}
\end{proposition}

\medskip
\noindent
Now let $\underline{\mu}$ be the partition obtained from $\mu$ by removing a box to the end of {\it every} row (or equivalently removing the first column), and let $\overline{\mu}$ be the partition obtained from $\mu$ by {\it adding} a box to every row. 

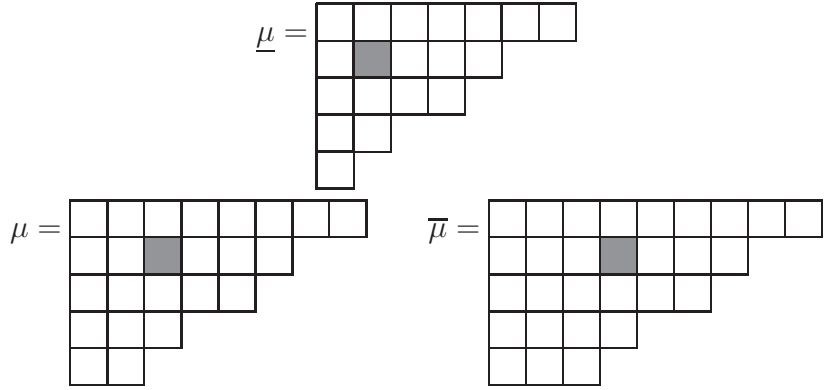
\begin{gather*}
\underline{\mu}=
\tableau{ \ &\ &\ &\ &\ &\ &\  \\ 
\ &\graycell &\ &\ &\  \\ 
\ &\ &\ &\   \\ 
\ &\  \\
\   }
\\
\mu = 
\tableau{ \ &\ &\ &\ &\ &\ &\ &\ \\ 
\ &\ &\graycell &\ &\ &\ \\ 
\ &\ &\ &\ &\  \\ 
\ &\ &\ \\
\ &\  }
\qquad
\overline\mu =
\tableau{ \ &\ &\ &\ &\ &\ &\ &\ &\  \\ 
\ &\ &\ &\graycell &\ &\ &\ \\ 
\ &\ &\ &\ &\ &\ \\ 
\ &\ &\ &\ \\
\ &\ &\ }
\end{gather*}


Although the highlighted box has the same arm and leg length in every partition, it has different co-ordinates with respect to the top left hand corner, giving us:
$$ H_{\underline{\mu}}(i-1,\mu_j) = H_\mu(i, \mu_j) = H_{\overline{\mu}}(i+1,\mu_j+1)$$
and similarly:
$$ \tilde{H}_{\underline{\mu}}(i-1,\mu_j) = \tilde{H}_\mu(i, \mu_j) = \tilde{H}_{\overline{\mu}}(i+1,\mu_j+1)$$

\medskip
\noindent
Observe further that $\overline{\mu}_-(\tilde{\alpha}) = \mu_+(\alpha)$ while $\underline{\mu}_+(\tilde{\alpha}) = \mu_-(\alpha)$. Using these facts we may rewrite the previous proposition as:

	
\begin{proposition}\label{fish}
We have:
\begin{align*}
\Phi(A_I,A_J)_{\epsilon q, t} &= \prod_{\substack{i \not\in \alpha \\ j \in \alpha \\ i < j}} \frac{\tilde{H}_{\mu_-(\alpha)}(i,\mu_j)}{\tilde{H}_\mu(i,\mu_j)}
\prod_{\substack{i \not\in \alpha \\ j \in \alpha \\ j < i}} \frac{H_{\mu_-(\alpha)}(j,\mu_i)}{H_\mu(j,\mu_i+1)}\\
\end{align*}
and:
\begin{align*}
\Phi(A_J,A_I)_{\epsilon q, t} &= \prod_{\substack{i \not\in \alpha \\ j \in \alpha \\ i < j}} \frac{H_{\mu_+(\alpha)}(i,\mu_j+1)}{H_\mu(i,\mu_j+1)} 
\prod_{\substack{i \not\in \alpha \\ j \in \alpha \\ j < i}} \frac{\tilde{H}_{\mu_+(\alpha)}(j,\mu_i+1)}{\tilde{H}_\mu(j,\mu_i)} \\
\end{align*}
\begin{proof}
This proposition reduces to the previous one after observing in the first case that:

$$ H_{\mu_-(\alpha)}(j,\mu_i) =H _{\underline{\mu}_+(\tilde{\alpha})}(j,\underline{\mu}_i+1) = H_{\mu_+(\tilde{\alpha})}(j,\mu_i+1)$$
and in the second case that:
$$ \tilde{H}_{\mu_+(\alpha)}(j,\mu_i+1) = \tilde{H}_{\overline{\mu}_-(\tilde{\alpha})}(j,\overline{\mu}_i) = \tilde{H}_{{\mu}_-(\tilde{\alpha})}(j,{\mu}_i)$$
\qedhere
\end{proof}
\end{proposition}

\subsubsection{Step three}

Using proposition \ref{fish} we may now give an alternative description for the expressions $R(\mu,\gamma)$ and $L(\mu,\gamma)$ arising in equation \eqref{toprove}. 

\begin{proposition} \label{RHS}
If $\gamma = \mu_-(\alpha)$ then:
$$ R(\mu,\gamma) =  w(1/t : A_J)_{\epsilon q, t} \,\, \Phi(A_I,A_J)_{\epsilon q, t} $$
\begin{proof}
By equation \eqref{one} we can write:
\begin{align*}
\prod_{s \in \tilde{R}_{\mu / \gamma} \cap C_{\mu / \gamma}} \frac{\tilde{H}_\gamma(s)}{\tilde{H}_\mu(s)} 
&= \prod_{\substack{ i \not\in \alpha \\ j \in \alpha \\ i < j }} 
\frac{\tilde{H}_{\mu_-(\alpha)}(i,\mu_j )}{\tilde{H}_{\mu}(i,\mu_j)} \\
&= \prod_{\substack{ i \not\in \alpha \\ j \in \alpha \\ i < j }} \Phi(a_i, a_j)_{\epsilon q, t} \\
\end{align*}
On the other hand, by shifting indices, we have:
\begin{align*}
\prod_{s \in R_{\mu / \gamma}} \frac{H_\gamma(s)}{H_\mu(s)} &=
\prod_{j \in \alpha} \left ( \prod_{k=1}^{\mu_j - 1} \frac{H_{\mu_-(\alpha)}(j,k)}{H_\mu(j,k)}\right ) \frac{1}{H_\mu(j,\mu_j)} \\
&= \prod_{j \in \alpha} \frac{1}{H_\mu(j,1)} \left ( \prod_{k=1}^{\mu_j-1} \frac{H_{\mu_-(\alpha)}(j,k)}{H_\mu(j,k+1)}\right ) \\
\end{align*}
Suppose that we take the diagram for $\gamma = \mu_-(\alpha)$ and place it on top of the diagram for $\mu$, but shifted to the right by one column. 

\[
\tableau{ \thickcell \mu &\ &\ &\ &\ &\ &\ &\ &\thickcell \gamma \\ 
\graycell \mu &\graycell &\ &\ &\graycell &\ \\ 
\thickcell \mu &\ &\ &\ &\thickcell \gamma \\ 
\graycell \mu &\graycell &\ \\
\thickcell \mu &\thickcell \gamma }
\]

\medskip
\noindent
For $j \in \alpha$, we have that $H_{\mu}(j,k+1) = H_{\mu_-(\alpha)}(j,k)$ unless $k=\mu_i$ for some $i \not\in \alpha$ with $i > j$. 

\medskip
\noindent
In the above diagram $\mu = (8,5,4,2,1)$ and $\alpha = \{2,4\}$. The grey boxes correspond to the terms in the expression which do not get cancelled.

\medskip
It follows then, by the first part of proposition \ref{fish}, that:
\begin{align*}
\prod_{j \in \alpha}  \prod_{k=1}^{\mu_j-1} \frac{H_{\mu_-(\alpha)}(j,k)}{H_\mu(j,k+1)}
&=  \prod_{\substack{j \in \alpha \\ i \not \in \alpha \\ i > j}} \frac{H_{\mu_-(\alpha)}(j,\mu_i)}{H_\mu(j,\mu_i+1)} \\
&= \prod_{\substack{j \in \alpha \\ i \not \in \alpha \\ i > j}} \Phi(a_i,a_j)_{\epsilon q, t}
\end{align*}

It remains to observe that in the first column we have:
$$\prod_{j \in \alpha} \frac{1}{H_\mu(j,1)} = \prod_{j \in \alpha} \frac{1 - a_j}{1 + t a_j / q} = w(1/t : A_J)_{\epsilon q,t} $$

\qedhere
\end{proof}
\end{proposition}

\bigskip\bigskip
\begin{proposition} \label{LHS}
If $\lambda = \mu_+(\alpha,p)$ then:
\begin{multline*}
L(\lambda,\mu) = \frac{(-q;t)_p}{(t;t)_p} \times \\
w(1/t : A_I)_{\epsilon q, t} V(1/t : t^{p+1} A_I)_{\epsilon q, t} W(1/t : t^p A_J)_{\epsilon q, t} \Phi(A_J,A_I)_{\epsilon q, t} 
\end{multline*}
\begin{proof}
The proof is similar to the previous one, only we must be slightly more careful with the first column. Begin by splitting, into two pieces, the expression:
$$ \prod_{s \in \tilde{R}_{\lambda / \mu} \cap C_{\lambda / \mu}} \frac{H_\lambda(s)}{H_\mu(s)} $$
The first piece is the contribution from the first column:
$$\prod_{\substack{i = 1 \\ i \not\in \alpha}}^n \frac{H_{\mu_+(\alpha,p)}(i,1)}{H_\mu(i,1)} $$
The second piece is everything else:
\begin{align*}
\prod_{\substack{ i \not\in \alpha \\ j \in \alpha \\ i < j }} \frac{H_{\mu_+(\alpha,p)}(i,\mu_j +1)}{H_\mu(i,\mu_j + 1)} 
&= \prod_{\substack{ i \not\in \alpha \\ j \in \alpha \\ i < j }} \frac{H_{\mu_+(\alpha)}(i,\mu_j +1)}{H_\mu(i,\mu_j + 1)} \\
&= \prod_{\substack{ i \not\in \alpha \\ j \in \alpha \\ i < j }} \Phi(a_j,a_i)_{\epsilon q, t} 
\end{align*}
Here we have made use of equation \eqref{three} to simplify the expression. We can ignore the extra ``tail'' of length $p$ because it only affects the first column.

\medskip
\noindent
Similarly, we split, into two pieces, the expression:
$$ \prod_{s \in R_{\lambda / \mu}} \frac{\tilde{H}_\lambda(s)}{\tilde{H}_\mu(s)} $$
Again the first piece is the contribution from the first column:
$$ \prod_{l=1}^p \tilde{H}_{\mu_+(\alpha)}(m+l,1) \prod_{j \in \alpha} \tilde{H}_{\mu_+(\alpha,p)}(j,1) $$
The second piece is everything else:
\begin{align*}
\prod_{j \in \alpha}  \prod_{k=1}^{\mu_j} \frac{\tilde{H}_{\mu_+(\alpha,p)}(j,k+1)}{\tilde{H}_\mu(j,k)} 
&= \prod_{j \in \alpha}  \prod_{k=1}^{\mu_j} \frac{\tilde{H}_{\mu_+(\alpha)}(j,k+1)}{\tilde{H}_\mu(j,k)} \\
\end{align*}
Again we can ignore the extra ``tail'' of length $p$ because it only affects the first column.

\bigskip\bigskip
\noindent
Suppose that this time we take the diagram for $\mu$ and place it on the top of the diagram for $\lambda = \mu_+(\alpha)$ but shifted to the right by one column. 

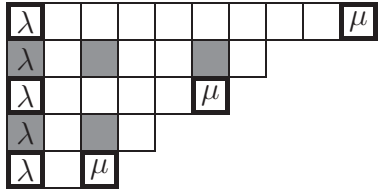
\[
\tableau{ \thickcell \lambda &\ &\ &\ &\ &\ &\ &\ &\ &\thickcell \mu \\ 
\graycell \lambda &\ &\graycell &\ &\ &\graycell &\ \\ 
\thickcell \lambda &\ &\ &\ &\ &\thickcell \mu \\ 
\graycell \lambda &\ &\graycell &\ \\
\thickcell \lambda &\ &\thickcell \mu \\
 }
\]

\medskip
\noindent
For $j \in \alpha$, we have that $H_{\mu_+(\alpha)}(j,k+1) = H_{\mu}(j,k)$ unless $k = 1$ or $k=\mu_i$ for some $i \not\in \alpha$ with $i > j$. 

\medskip
\noindent
In the above diagram $\mu = (8,5,4,2,1)$ and $\alpha = \{2,4\}$. Again the grey boxes correspond to the terms which do not get cancelled. 

\medskip
We have, by the second part of proposition \ref{fish}, that:
\begin{align*}
\prod_{j \in \alpha} \prod_{k=1}^{\mu_j} \frac{\tilde{H}_{\mu_+(\alpha)}(j,k+1)}{\tilde{H}_\mu(j,k)}
&= \prod_{\substack{j \in \alpha \\ i \not \in \alpha \\ i > j}} \frac{\tilde{H}_{\mu_+(\alpha)}(j,\mu_i+1)}{\tilde{H}_\mu(j,\mu_i)} \\
&=  \prod_{\substack{j \in \alpha \\ i \not \in \alpha \\ i > j}} \Phi(a_j, a_i)_{\epsilon q, t} \\
\end{align*}
It remains to find suitable expressions for the four remaining factors coming from the first column:
\begin{align*}
 \prod_{l=1}^p \tilde{H}_{\mu_+(\alpha,p)}(m+l,1) &= \Omega( (t - \epsilon q) (1 + t + \cdots + t^{p-1})) = \frac{(-q;t)_p}{(t;t)_p} \\
 \prod_{j \in \alpha} \tilde{H}_{\mu_+(\alpha,p)}(j,1) &= \prod_{j \in \alpha} \frac{1 + qt^p a_j}{1 - t^{p+1} a_j} = W(1/t : t^{p} A_J)_{\epsilon q, t} \\
\prod_{\substack{i = 1 \\ i \not\in \alpha}}^m H_{\mu_+(\alpha,p)}(i,1) &= \prod_{\substack{i = 1 \\ i \not\in \alpha}}^m \frac{1 + t/q \,\, a_i \, t^p}{1 - a_i \, t^{p}} = V(1/t : t^{p-1} \,\, A_I)_{\epsilon q, t}\\
 \prod_{\substack{i = 1 \\ i \not\in \alpha}}^m \frac{1}{H_\mu(i,1)} &= \prod_{\substack{i = 1 \\ i \not\in \alpha}}^m  \frac{1 - a_i}{1 + t/q \,\, a_i } = w(1/t : A_I)_{\epsilon q, t}
\end{align*}

\qedhere
\end{proof}
\end{proposition}

\bigskip\bigskip
\begin{proposition}
The Kawanaka identity is equivalent to the following:
\begin{multline*}
\!\!\!\!\sum_{s = 0}^k \frac{(q;t)_s}{ (t;t)_s} \sum_{\substack{I \cup J = [m] \\ |J| = k-s}} \Big(w(1/t : A_I) V(1/t : t^{s-1} A_I) W(1/t : t^{s} A_J) \Phi(A_J,A_I) 
\\- w(1/t :  A_J) \Phi(A_I,A_J) \Big) = 0
\end{multline*}

\begin{proof}
Equation \eqref{toprove} may be rewritten in the form:
$$ \sum_{\lambda \in \tilde{U}_k(\mu)} L(\lambda, \mu) - \sum_{s = 0}^k \frac{(-q;t)_s}{(t;t)_s} \sum_{\gamma \in \tilde{D}_{k-s}(\mu)} R(\mu, \gamma) = 0 $$
If $\mu$ is a partition with distinct parts, then, as remarked at the end of the proof of proposition \ref{schur}, we may replace the sum over vertical strips with a sum over subsets of the rows, to obtain the equivalent expression:
$$ \sum_{s = 0}^k \sum_{|\alpha| = k-s} \left( L(\mu_+(\alpha,s), \mu) -  \frac{(-q;t)_s}{(t;t)_s} R(\mu, \mu_-(\alpha)) \right) = 0 $$

\medskip
The result now follows, in this special case, by applying propositions \ref{RHS} and \ref{LHS} and then replacing $\epsilon q$ with $q$ (formal negative).

\medskip
\noindent
If $\mu$ contains a repeated part, and $\alpha$ is such that $\mu_-(\alpha)$ is not a partition then we have $a_{i-1} = ta_i$ for some $j=i-1 \in \alpha$ and $i \not \in \alpha$ and so the expression:
\begin{align*}
\Phi(A_I,A_J)
&= \frac{(a_i + a_j t/q)(a_i - a_j/t)}{(a_i - a_j)(a_i + a_j/q)}
\end{align*}
vanishes.

\medskip
\noindent
Similarly if $\mu$ contains a repeated part, and $\alpha$ is such that $\mu_+(\alpha,p)$ is not a partition then we have $1/ta_i = a_{i+1}$ for some $i \not \in \alpha$ and $j = i+1 \in \alpha$ and the expression: 
\begin{align*}
\Phi(A_J,A_I) 
&= \prod_{\substack{i \not \in \alpha \\ j \in \alpha}}\frac{(a_j + a_i t/q)(a_j - a_i/t)}{(a_j - a_i)(a_j + a_i/q)}
\end{align*}
vanishes.

\medskip
\noindent
Thus, even when $\mu$ contains a repeated part, we may still replace the sum over vertical strips with a sum over subsets of the rows, since the terms which do not correspond to partitions all vanish. The result follows.

\end{proof}
\end{proposition}

The above identity is none other than a special case of proposition \ref{final} with $z = 1/t$, which completes the proof of the Kawanaka identity:

$$ \sum_\lambda  \prod_{s \in \lambda}  \left ( \frac{1+q^{a_\lambda(s)}t^{l_\lambda(s)+1}}{1-q^{a_\lambda(s)+1}t^{l_\lambda(s)}} \right )
P_\lambda(X;q^2,t^2) = \prod_{i \geq 1} \frac{(-tx_i;q)_\infty}{(x_i;q)_\infty} \prod_{i<j} \frac{(t^2x_ix_j;q^2)_\infty}{(x_ix_j;q^2)_\infty} $$

\nocite{SPECIES}
\nocite{COLROT}	
\bibliographystyle{amsplainhyper} 
\bibliography{references}

\end{document}